%% file: main.tex
\documentclass[12pt]{article}

\input{header}
\begin{document}

\maketitle

\begin{abstract}
  We consider solving stochastic programs over an infinite horizon. 
  By leveraging the stationarity of the problem, we develop a novel continually-exploring infinite-horizon explorative dual dynamic programming (CE-Inf-EDDP) algorithm.
  CE-Inf-EDDP builds upon the existing explorative dual dynamic programming, designed for the finite-horizon problem, by specializing it for the infinite-horizon, stationary case. 
  \revone{By incorporating cut sharing, frequent cutting-plane model updates, and a new adaptive search point selection strategy, CE-Inf-EDDP provides state-of-the-art iteration complexity while offering encouraging numerical performance.}
  In the newsvendor and hydrothermal planning problem, CE-Inf-EDDP can reduce the runtime of each iteration by up to one to two orders of magnitude compared to prior methods while maintaining similar solution quality. 
  As a result, the final solution quality and its guarantees can be much improved over the same runtime.
  For example, in the hydrothermal planning problem, CE-Inf-EDDP attains about an order of magnitude improvement in the relative optimality gap compared to existing methods.
\end{abstract} 

\vspace{0.5em}


\vspace{1em}
\mscclass{49M37, 90C15, 90C25, 90C39, 93A13}

\section{Introduction}
Making decisions under uncertainty over a multi-period horizon is a fundamental problem that has been extensively studied through the lens of Markov decision processes, stochastic optimal control, and stochastic programming~\cite{lan2024numerical}.
Markov decision processes can be efficiently solved when the number of states and actions is relatively small, but extensions to large or continuous state and action spaces require function approximation~\cite{powell2007approximate}, which may yield sub-optimal solutions.
Some specially structured stochastic optimal control problems (e.g., linear dynamics with quadratic costs~\cite{caines2018linear}) admit computable closed-form solutions, but incorporating general constraints or costs is challenging.
In contrast, multi-stage stochastic programs can potentially handle both continuous state and action spaces as well as general constraints and cost functions.
However, existing results typically focus on the finite-horizon case.
In reality, decisions must be made with consideration over an infinite horizon, such as in energy operations~\cite{pereira1991multi,cruise2019control}, production planning~\cite{mcclain1977horizon,van2005comparison}, and finance~\cite{jurek2011optimal,brown2011dynamic}.

Therefore, we focus on an important subclass of \textit{stationary infinite-horizon stochastic programs}, which is the multi-stage stochastic program with discounted costs and the number of stages $T \to \infty$. 
The cost, (parametric) feasible set, and distribution over random variables are identical across all stages.
\revone{We assume the randomness, captured by random variable $\xi$, is stage-wise independent.
This assumption appears common in the literature~\cite{shapiro2021lectures,philpott2008convergence,shapiro2011analysis} to ensure computational tractability.
This assumption, together with the stationarity of the problem, implies that each stage's random vector is independent and identically distributed (i.i.d.).}
Following the dynamic programming formulation, i.e., Wald-Bellman equations~\cite[Section 3]{shapiro2020periodical}, the problem can be stated succinctly as
\begin{align} \label{eq:qgp_a1}
  \min_{x^1 \in X(x^0,\xi^1)} \{h(x^1,\B c^1) + \lambda \mathcal{V}(x^1) \},
\end{align}
where $x^0$ and $\xi^1 := \{\B A^1,\B B^1,\B Q^1, \B b^1, \B c^1, \B r^1\}$ are fixed vectors, 
$x^1$ is a decision variable referred to as the \textit{here-and-now} decision, $\lambda \in (0,1)$ is a discount factor, and $\mathcal V$ is the so-called \textit{value function} (defined in the next paragraph).
The parametric feasible set is 
\begin{talign} \label{eq:qgp_a01}
  X(u,\xi) 
  &= 
  \{ x \in \bar{X} \subseteq \mathbb{R}^n: \B{A}x - \B{B}u - \B b  \in \mathcal{K}, \phi(x, \B r) \leq \B{Q} u\},
\end{talign}
where $\xi := \{\B A, \B B, \B Q, \B b, \B c, \B r\}$ is a random vector.
The cost function $h(\cdot, \B c): \mathbb{R}^n \to \mathbb{R} $ is a closed, proper convex function parameterized by vector $\B c$, $\bar{X} \subseteq \mathbb{R}^n$ is a closed convex set, and $\mathcal{K}$ is a regular cone.
Here, $\B{A} : \mathbb{R}^{n} \to \mathbb{R}^{m}$, $\B{B} : \mathbb{R}^n \to \mathbb{R}^m$, and $\B{Q} : \mathbb{R}^n \to \mathbb{R}^{m_\phi}$ are linear mappings for some integers $m \geq 0$ and $m_\phi \geq 0$, while $\phi_i(\cdot, \B r) : \mathbb{R}^n \to \mathbb{R}$ is a closed convex functional constraint for all $i = 1,\ldots,m_\phi$ parameterized by $\B r$.

The value function $\mathcal V(x)$ quantifies the future costs of a here-and-now decision $x$, which satisfies the so-called \textit{Bellman optimality equations},
\begin{align} \label{eq:qgp_a2}
  \mathcal V(x) &:= \mathbb{E}_{\xi}[\tilde{\mathcal{V}}(x,\xi)] \\
  \tilde{\mathcal{V}}(u,\xi) &:= \min_{x \in X(u,\xi)} \{ h(x,\B c) + \lambda \mathcal V(x) \}, \nonumber
\end{align}
where the equations should hold over the effective domain of $\mathcal V$.
\revone{Our stationary, infinite-horizon formulation can be viewed as a special case of the policy graphs of~\cite{dowson2020policy} with a cyclic graph or of the periodical problems of~\cite{shapiro2020periodical} with periodicity $\mathrm{m}=1$. 
Our intention is to show that a careful specialization of existing methods towards the aforementioned problem can achieve a better theoretical complexity and deliver promising numerical performance over the general-case algorithms.
}

This paper presents a new class of dual dynamic programming (DDP) methods for solving~\eqref{eq:qgp_a1}.
Our methods essentially tailor the explorative dual dynamic programming (EDDP)~\cite{lan2022complexity} for solving finite-horizon problems.
Recall that EDDP is a deterministic counterpart of the celebrated stochastic dual dynamic programming (SDDP) as first introduced by Pereira and Pinto~\cite{pereira1991multi}, the latter of which is based on the nested cutting-plane method by Birge~\cite{birge1985decomposition}.
In such nested cutting-plane methods for solving convex multi-stage programs, one builds approximations of the value function via a cutting-plane model.
This cutting-plane model is utilized in the \textit{forward phase}, where a series of subproblems are solved to find the most promising points, or \textit{trial points}, in the feasible region to explore.
Then in the \textit{backward phase}, the trial points are used to improve the cutting-plane model.
Then the algorithm repeats these two phases until termination.
Note that EDDP can be applied directly towards stationary, infinite-horizon problems by properly setting a user-defined horizon length $T$.
In the context of infinite-horizon problems, $T$ will be referred to as the \textit{effective planning horizon}, since it provides a trade-off between accuracy and complexity for a finite-horizon model to approximate an infinite-horizon one.

Our main contribution is the design, analysis, and implementation of a new DDP method called continually-exploring infinite-horizon EDDP, or CE-Inf-EDDP.
\revone{First, CE-Inf-EDDP reduces both the computational and the memory costs of each iteration by $T$ compared to EDDP.
This is achieved in two ways.
For the first way, CE-Inf-EDDP reduces the number of cutting-plane models from $T$ to $1$ through a technique called cut sharing, which has also appeared in prior works~\cite{dowson2020policy,shapiro2020periodical,nannicini2021benders}.
For the second way, the method updates this cutting-plane model more frequently by reducing the number of steps in the forward phase from $T$ to $1$, allowing the model to be updated before visiting the next stage.
While the more frequent model updating is not new in DDP methods (see~\cite{chen1999convergent} for an example), we show that combining it with cut sharing yields significant numerical improvements.}
\revone{Second, CE-Inf-EDDP incorporates a new search point selection that explores the feasible region longer.
We show this new selection improves the iteration complexity of CE-Inf-EDDP by a factor of $T$.
Third, we provide new variants of CE-Inf-EDDP with different search point selection rules. 
In particular, our alternative search point selection rules include one with an upper-bound model and one with random selection, which is similar to SDDP.
We also study the convergence of these algorithms under inexact solutions. 
}
Lastly, we provide the first implementation of EDDP-like algorithms (with open-source code\footnote{\url{https://github.com/jucaleb4/hddp}}). Our results are promising. 
For example, our methods offer multiple orders of magnitude improvement in the optimality gap compared to existing DDP-type methods.

Although some recent developments in EDDP-like algorithms also improve the iteration complexity by $T$~\cite{zhang2020distributionally}, the improvements are more complicated because they require both a lower-bound model (cutting-plane model) and an upper-bound model, whereas CE-Inf-EDDP only needs a cutting-plane model.
Finally, we mention that our proposed methods are slightly more general than existing ones since our algorithm design and analysis can handle nonlinear, convex costs and constraints.
In contrast, existing non-asymptotic convergence analyses are often limited to the linear case~\cite{shapiro2023dual,baucke2017deterministic,guigues2020inexact}.

The rest of the paper is organized as follows.
In the next Section~\ref{sec:sec2}, we introduce our infinite-horizon dual dynamic programming methods, including CE-Inf-EDDP.
In Section~\ref{sec:sec4}, we consider several practical extensions or modifications to CE-Inf-EDDP to monitor or improve convergence in practice.
Finally, in Section~\ref{sec:experiments}, we present preliminary numerical experiments.

\section{Infinite-horizon explorative dual dynamic programming} 
\label{sec:sec2}

Throughout this paper, we make the following two assumptions regarding the main problem of interest,~\eqref{eq:qgp_a1}.
The first assumption pertains to the effective feasible region.
Let us define the feasible region for the $t$-th stage, recursively defined by
\begin{talign*}
  \mathcal{X}^t 
  =
  \begin{cases}
    X(x^0,\xi^1) & : t = 0 \\
    \cup_{x \in {\mathcal{X}}^{t-1}} \cup_{\xi \in \Theta} X(x, \xi) & : t \geq 1
  \end{cases},
\end{talign*}
where $\Theta$ is the space of all random vectors $\xi = \{A,B,Q,b,c,r\}$.
We also define the union across all stages, $\mathcal X := \cup_{t \geq 0} \mathcal{X}^t$.
\begin{assumption}  \label{asmp:eddp_asmp_domain}
  There exists a constant $D < +\infty$ such that,
  \[
    \|x - x'\|_\infty \leq D, \ \forall x,x' \in \mathcal X.
  \]
\end{assumption} 
Note that this assumption can be easily satisfied if the set $\bar{X}$ from the parametric feasible set~\eqref{eq:qgp_a01} is bounded.
Moreover, this assumption seems standard in both the analysis of the sample average approximation problem~\cite{shapiro2021lectures} and cutting-plane type methods~\cite{lan2022complexity,nesterov2018lectures}.

To assist our regularity conditions later, let $\mathrm{Aff}(\mathcal X)$ be the affine hull of the corresponding set and let the $\bar{\epsilon}$-\revtwo{hypercube} surrounding it be $\mathcal B(\mathcal X, \bar \epsilon) := \{y \in \mathrm{Aff}(\mathcal X) : \|y\|_\infty \leq \bar{\epsilon}\}$ for some $\bar{\epsilon} > 0$. 
Let
\begin{align}\label{eq:qgp_a9}
  \mathcal X(\bar \epsilon) = \mathcal X + \mathcal{B}(\mathcal X, \bar{\epsilon}).
\end{align}
\revone{The space $\mX(\bar \epsilon)$ is only used to establish certain regularity conditions, such as the following assumption regarding the relative interior (rint), needed to establish convergence for our algorithm. 
Therefore, $\mX(\bar \epsilon)$ will not be directly used in the algorithm design.}

\begin{assumption}  \label{asmp:eddp_asmp_relint}
  There exists $\bar \epsilon \in (0,+\infty)$ s.t. $\mathrm{rint}(X(u,\xi) ) \ne \emptyset, \ \forall u \in \mathcal X(\bar \epsilon), \ \xi \in \Theta$.
\end{assumption} 
This assumption implies every ``min'' in~\eqref{eq:qgp_a2} is feasible.
Such an assumption is related to the relatively complete recourse, and it appears to be commonly made in the multi-stage stochastic programming literature~\cite{shapiro2021lectures,shapiro2011analysis,dowson2020policy,philpott2008convergence}.
It is also closely related to the extended relatively complete recourse assumption~\cite{girardeau2015convergence}.
The requirement of a nonempty relative interior implies the Slater condition, ensuring the existence of an optimal dual solution~\cite{lan2022complexity} to construct subgradients for a cutting-plane model.
If the conic constraints in~\eqref{eq:qgp_a1} are polyhedrally representable, e.g., $\mathcal{K} = \{\mathbf{0}\}$ or $\mathbb{R}^n_+ = \{x \in \mathbb{R}^n : x_i \geq 0\}$, then by the relaxed Slater conditions~\cite{ben2001lectures} we only require $X(u, \xi) \ne \emptyset$. 

Still, solving~\eqref{eq:qgp_a1} is challenging.
One reason is that $\mathcal V$ is defined as an infinite-nested sequence of minimizing stochastic programs involving expectations (i.e., multi-dimensional integrals), which are not known a priori. 
We first resolve the intractability from the expectation by approximating it via a sample average approximation, or SAA~\cite{kleywegt2002sample}.
Let us draw $N$ i.i.d.~realizations of the random variable $\xi$, denoted by $\{ \tilde{\xi}_i\}_{i=1}^N$, where
\begin{talign} \label{eq:saa_xi}
    \tilde{\xi}_i 
    := 
    \{\tilde{\B{A}}_i, \tilde{\B{B}}_i, \tilde{\B{Q}}_i, \tilde{\B{b}}_i, \tilde{\B{c}}_i, \tilde{\B{r}}_i\}.
\end{talign}
The resulting SAA problem is
\begin{talign}  
\label{eq:qgp_a3}
  \min_{x^1 \in X(x^0,\xi^1)} \{ F(x) := h(x^1,\B c^1) + \lambda V(x^1)\},
\end{talign}
where
\begin{talign} \label{eq:qgp_a3b}
   V(x) &:= \frac{1}{N} \sum_{i=1}^N v_{i}(x) \\
  v_{i}(u) &:= \min_{x \in X(u, \tilde{\xi}_{i})} 
                \left \{ F_i(x) := h(x,\tilde{\B{c}}_i) + \lambda V(x) \right \}. \nonumber
\end{talign}

We make a few remarks.
First, the definition of the value function $V$ differs from $\mathcal V$ in~\eqref{eq:qgp_a2} in that the former uses an empirical average over the random variable $\xi$, while the latter uses the true expectation.
Second, since the realized random vectors $\tilde{\xi}_i$ are drawn from the distribution with respect to (w.r.t.)~the random variable $\xi$, then clearly Assumption~\ref{asmp:eddp_asmp_relint} still holds w.r.t.~$\tilde{\xi}_i$ in place of $\xi$.
Third, we observe an important Lipschitz property of the sampled cost function.
With some abuse of notation, we write $\tilde{c}_0 = \B c^1$ from~\eqref{eq:qgp_a1} to simplify notation.
\begin{corollary} \label{corr:h_minimizes}
    For any fixed $\B c \in \xi \in \Theta$ (such as $\tilde{\B c}_i$ and $c^1$) and any $\bar \epsilon \in (0,+\infty)$, the cost function $h(\cdot, \B c)$ is bounded over $\mathcal{X}(\bar \epsilon)$, and in particular,
    \begin{talign*}
      -\infty < \underline{h} := \min_{i \in\{0,\ldots,N\},x \in \mathcal{X}(\bar \epsilon)} h(x, \tilde{\B c}_i) 
      \leq
      \max_{i \in \{0,\ldots,N\},x \in \mathcal{X}(\bar \epsilon)} h(x, \tilde{\B c}_i)
      =: \overline{h}
      <+\infty.
    \end{talign*}
    Also, it is Lipschitz continuous over $\revtwo{\mathcal{X}}$ w.r.t.~the norm $\|\cdot\|_\infty$ with constant $M_h \leq \revone{\frac{\overline{h} - \underline{h}}{\bar \epsilon}}$, or
    \begin{talign*}
      \vert h(x, \B c) - h(x', \B c) \vert \leq M_h \|x'-x\|_\infty, \ \forall x \in \revtwo{\mathcal X}.
    \end{talign*}
\end{corollary}
\begin{proof}
    Any proper Lipschitz function over a compact domain is bounded~\cite[Theorem 10.6]{rockafellar1972convex}.
    Since $h(\cdot,\B c)$ is proper convex over $\mathbb{R}^n$ and $\mathcal{X} \subset \mathrm{rint}(\mathcal{X}(\bar \epsilon))$ is compact (Assumption~\ref{asmp:eddp_asmp_domain}), the function is Lipschitz continuous over $\mathcal{X}$ w.r.t.~the norm $\|\cdot\|_\infty$, and the bound on $M_h$ can be derived from~\cite[Theorem 10.4]{rockafellar1972convex}.
\end{proof}
The corollary states the Lipschitz continuity w.r.t.~the norm $\|\cdot\|_\infty$.
For the rest of the paper, we assume Lipschitz continuity is defined w.r.t.~the norm $\|\cdot\|_\infty$.

Let us briefly discuss how to obtain bounds $\underline{h}$ and $\overline{h}$, since we will need them to design algorithms later.
Since $h(\cdot, \tilde{\B c}_i)$ is convex, $h$ can be minimized over the known set $(\bar X + \mathcal{B}(\bar X, \bar \epsilon)) \supseteq \mathcal{X}(\bar \epsilon)$, where recall $\bar X$ from~\eqref{eq:qgp_a01}. 
On the other hand, computing $\overline{h}$ may be difficult to obtain, since maximizing a convex function is NP-hard in general.
Using the Lipschitz constant $M_h$, one can obtain the crude bound $\overline{h} = \underline{h} + M_h D$.
To end, we state a useful property for the value function $V$.
\begin{lemma} \label{lem:eddp_asmp_smooth1}
    The value function $V$ is convex over $\mathcal{X}(\bar \epsilon)$.
    Furthermore, $V$ is Lipschitz continuous over $\mathcal X$ with constant $M_V \leq \revone{\frac{\overline{h} - \underline{h}}{(1-\lambda) \bar \epsilon}}$, or
    \begin{talign*}
        \vert V(x') - V(x) \vert
        \leq
        M_V\|x'-x\|_\infty,
        \ \forall x',x \in \mathcal X.
    \end{talign*}
\end{lemma}
\begin{proof}
    Convexity of $V$ can be shown similarly to~\cite[Proposition 1]{grinold1977finite}.
    It remains to show $V(x)$ is Lipschitz continuous over $\mathcal X \subseteq \mathrm{rint}(\mathcal X (\bar \epsilon))$.
    We can apply the same argument as Corollary~\ref{corr:h_minimizes}, where it suffices to show $V(\cdot)$ is bounded, and hence proper, over $\mathcal X(\bar \epsilon)$.
    By Assumption~\ref{asmp:eddp_asmp_relint}, $V$ is well-defined over $\mathcal{X}(\bar \epsilon)$.
    In particular, Corollary~\ref{corr:h_minimizes} ensures that the value function is absolutely bounded with $V(x) \in [\frac{\underline{h}}{1-\lambda}, \frac{\overline{h}}{1-\lambda}]$ for any $x \in \mathcal{X}(\bar \epsilon)$, as desired. 
\end{proof}
With these properties in hand as well as Assumption~\ref{asmp:eddp_asmp_domain},
one can show for a sufficiently large sample size $N$ that $V$ approximates the true value function $\mathcal V$ well~\cite[Theorem 5.35]{shapiro2021lectures}.
\revone{In particular, although $V$ will be a biased estimate of $\mathcal{V}$, optimizing w.r.t.~$V$ will output solutions that are nearly optimal for the unknown $\mathcal{V}$, which can be shown using SAA-type arguments~\cite[Lemma 5.18]{shapiro2021lectures}.}
Thus, for the remainder of our paper, we focus on the SAA problem~\eqref{eq:qgp_a3}.
To deal with the infinite nesting of minimization problems, we next present a stage-wise decomposition method.

\subsection{Background on EDDP} \label{sec:background_eddp}
\revone{
We start by reviewing EDDP from~\cite{lan2022complexity}, which our main algorithm builds upon.
Similar to the nested Bender's decomposition and stochastic dual dynamic programming, EDDP aims to solve a finite-horizon problem of length $T$.
EDDP is an iterative algorithm where each iteration $k$ involves a forward phase followed by a backward phase.
In the forward phase, feasible solutions are generated using the cutting-plane model built in previous iterations, where we start from the first stage and move forward until the last ($T$-th) stage.
Without making specific mention of the stage, let us denote a feasible solution at iteration $k$ as the \textit{search point} $x^k$.
We also write the cutting-plane model at iteration $k$ as $\uV^k : \mathcal{X} \to \mathbb{R}$.
Later, in Algorithm~\ref{alg:alg2}, we will provide a precise definition, and we will verify $\uV^k$ lower bounds $V(\cdot)$. 
After the forward phase, we run the backward phase, where the cutting-plane model is updated after each stage, starting from the last stage and going backwards until the first stage.}

\revone{
EDDP introduced a key concept called saturated points, which is also essential in our analysis.
Let us denote a sequence of accuracy tolerances $\{\epsilon_t\}_{t=0}^{T}$, which is assumed to be a non-decreasing sequence so that $\epsilon_0$ is the smallest value.
\begin{definition} \label{def:dqgp_a0}
  A search point $x^{k-1}$ is $\epsilon_t$-saturated at iteration $k$ when $V(x^{k-1}) - \uV^k(x^{k-1}) \leq \epsilon_t$.
\end{definition}
The saturation level provides an upper bound on the error between the cutting-plane model and the true value function at the search point $x^k$ at iteration $k$.
Let us denote $S^k_t$ as the set of search points up to iteration $k$ that are $\epsilon_t$-saturated.
Clearly, as EDDP finds more saturated points, or equivalently as $S^k_t$ becomes larger, $\uV^k$ becomes a provably better approximation of $V$.
The convergence analysis of EDDP shows that new saturated points will be added to $S^k_t$ for different $t$ as the number of iterations $k$ increases.
}

\revone{
However, EDDP is designed for the general finite-horizon case, where each stage can have a different cost and a parameterized feasible region.
Therefore, a direct application of EDDP to the stationary, infinite-horizon problem cannot properly exploit the fact that the value function $V(\cdot)$ is the same across all stages in~\eqref{eq:qgp_a3}.
In addition, EDDP requires storing the set of all saturated points $S^k_t$ to select new search points, which incurs both computational (in terms of an additional inner loop over forward and backward steps) and memory costs that scale linearly with both the iterations $k$ and stages $T$.
For infinite-horizon problems, $T$ must be set on the order of $O(\frac{1}{1-\lambda})$, which is large for a large discount factor $\lambda$.
We next introduce Inf-EDDP, where we specialize EDDP towards the infinite-horizon, stationary case with improvements in both the computational and memory costs.
}

\subsection{The basic Inf-EDDP} \label{sec:preliminary_eddp}
We now present our new algorithm, the infinite-horizon explorative dual dynamic programming, or Inf-EDDP for short, in Algorithm~\ref{alg:alg2}.
Similar to EDDP, Inf-EDDP takes as inputs a discretization level $\epsilon > 0$ and a horizon length $T$, which we refer to as the \textit{effective planning horizon}.
The choice of $\epsilon$ and $T$ will primarily affect the final optimality gap.
Inf-EDDP also has a cutting-plane model, denoted by $\uV^k$ at iteration $k$, which is iteratively updated and used to find new search points.

\begin{algorithm}[h!]
  \caption{}
  \label{alg:alg2}
  \begin{algorithmic}[1]
  \Procedure{Inf-EDDP}{$\epsilon, T$}
      \State{Initialize $\uV^0(x) = (1-\lambda)^{-1} \underline{h}, \ \forall x \in \mathbb{R}^n$} \label{line:line0} 
      \State{Set $\mathcal{S}_\epsilon(x) = T, \ \forall x \in \mXeps$ } \label{line:line1} \Comment{$\mXeps$ is the $\epsilon$-net of $\mX$}
      \For{$k=1,2,\ldots$}
          \State{Solve
          \vspace{-0.5em}
          \begin{talign} \label{eq:policy}
            x^k_0 \in \mathrm{argmin}_{x \in X(x^0,\xi^1)} \{ \uF^k(x) := h(x, \B c^1) + \lambda \uV^{k-1}(x) \}.
          \end{talign}} \label{line:line_fwd}
          \If{$\Seps(\Geps(x^k_0)) \leq 1$} \label{line:line6}
          \Comment{$\Geps(\cdot)$ is defined in~\eqref{eq:gap_eps}}
              \State{Terminate and return $x^k_0$} \label{line:line7} 
          \EndIf
          \For{$i=1,2,\ldots,N$}
              \State{
                Let $x^k_i$ and $(\tv^{k}_i)'(x^{k-1})$ be the primal sol'n and subgradient at $x^{k-1}$ for
                  \begin{align} \label{eq:eq5} 
                  \tv^{k}_i(x^{k-1}) 
                  = 
                  \min_{x \in X(x^{k-1}, \tilde{\xi}_i)}
                    \{ 
                    \uF^{k}_i(x) := h(x,\tilde{c}_i) + \lambda \uV^{k-1}(x) 
                    \}. 
                 \end{align}
              } \label{line:line2}
          \EndFor
          \State Update 
          \begin{talign} \label{eq:uV_max}
            \uV^k(x) = \max \big \{ \uV^{k-1}(x), \frac{1}{N} 
            \sum_{i=1}^N \left [ \tv_i^k(x^{k-1}) + \langle  (\tv^{k}_i)'(x^{k-1}), x-x^{k-1} \rangle \right ] 
            \big\}
          \end{talign} \label{line:line3b}
          \State{Set $z^k = \mathrm{argmax}_{\{x^k_i\}_{i=1}^N} \Seps(\Geps(x^k_i))$} \label{line:z_def}
          \State{Set} \label{line:line8}
          \begin{align}  \label{eq:next_pt}
            x^k = \begin{cases}
              \revone{x^k_0}  & : k \text{ mod } T = \revone{1} \\ 
              z^k & : \text{o/w}
            \end{cases}.
          \end{align}
          \State{Update $\Seps(\Geps(x^{k-1})) =  \min( \Seps(\Geps(x^{k-1})), \Seps(\Geps(z^k))-1)$} \label{line:line5}
      \EndFor
  \EndProcedure
  \end{algorithmic}
  \end{algorithm}

\revone{There are two key differences between EDDP and Inf-EDDP. 
First, Inf-EDDP uses one cutting-plane model $\uV^k(\cdot)$ for all stages, also known as cut sharing~\cite{dowson2020policy,shapiro2020periodical,nannicini2021benders}, while EDDP has one for each of the $T$ stages.
Second, instead of a forward and backward phase in EDDP, Inf-EDDP has a simpler forward and backward step, which we define in the next paragraph.
The forward (and similarly backward) phase from EDDP can be viewed as a sequence of $T$ forward steps. 
Therefore, the number of forward and backward steps in Inf-EDDP is reduced by a factor of $T$.
These two innovations yield more frequent updates to a single cutting-plane model, improving its overall accuracy.
}

\revone{
The \textit{forward step} at iteration $k$ consists of Line~\ref{line:line_fwd} to Line~\ref{line:line2}, which solves for a trial search point $x^k_i$ and dual solution (to compute the subgradient, which exists under Assumption~\ref{asmp:eddp_asmp_relint}) for sub-problems $i=0,\ldots,N$.
Sub-problem $i=0$ corresponds to the here-and-now problem~\eqref{eq:qgp_a3}, while the latter $N$ sub-problems correspond to scenario $i$ of the SAA problem~\eqref{eq:qgp_a3b} with a feasible region parameterized by the last search point, $x^{k-1}$.
In contrast, the forward phase in EDDP solves a single here-and-now sub-problem in the first stage and solves $N$ sub-problems for each stage $t=2,\ldots,T$.
The forward step also involves a stopping criterion in \revtwo{Line~\ref{line:line6}}, which we will discuss after introducing the backward step.}

\revone{
The \textit{backward step} at iteration $k$ consists of Line~\ref{line:line3b} to Line~\ref{line:line5}.
The cutting-plane model is updated in Line~\ref{line:line3b} using the latter $N$ trial search points and their associated subgradients (and it is initialized in Line~\ref{line:line0}). 
The next search point $x^k$ is selected from the initial search point $x^k_0$ every $T$ iterations and otherwise as an auxiliary search point $z^k$ in Line~\ref{line:line8}, and the stopping criterion is updated in Line~\ref{line:line5}.
The sequence $\{z^k\}_k$, where $z^k$ is explained later below, is the trajectory of search points found by repeatedly solving the sub-problems with the cutting-plane model.
Then, periodically every $T$ iterations, we reset the trajectory to $x^k_0$, which mimics the end of a finite-horizon problem of length $T$.
Similar to the forward step, the backward step simplifies the backward phase from EDDP by removing the pass over the update of $T$ cutting-plane models and $T$ stopping criteria.
}

\revone{
We will now discuss a stopping criterion, initialized in Line~\ref{line:line1}, queried in Line~\ref{line:line6},~\ref{line:z_def},~and~\ref{line:line5}, and also updated in Line~\ref{line:line5}. 
Define the \textit{saturation table},
\begin{align} \label{eq:S_def}
  \Seps: \mXeps \to \mathbb{Z},
\end{align}
where $\mXeps$ is the $\epsilon$-net of $\mX$ (w.r.t.~the distance metric, $d(\cdot,\cdot) = \|\cdot-\cdot\|_\infty$) and $\mathbb{Z}$ is the set of all integers.
We also define the projection operator w.r.t.~$\mXeps$ as the closest point in the $\epsilon$-net,
\begin{talign} \label{eq:gap_eps}
  \Geps(x) := \mathrm{argmin}_{s \in \mXeps} \|s-x\|_\infty,
\end{talign}
where ties can be broken arbitrarily but consistently.
The mapping $\Geps(\cdot)$ extends the domain of $\Seps$ from the finite set $\mXeps$ to the uncountable feasible region $\mX$.
Later, Proposition~\ref{prop:eddp_next_sat} will show that at iteration $k$, $\Seps(\Geps(x))$ upper bounds the lowest saturation level of a search point $x$ via
\begin{talign} \label{eq:Seps_ub}
  \Seps(\Geps(x)) \geq \min_{t} \{ t ~:~ V(x) - \uV^{k-1}(x) \leq \epsilon_{t-1}/\lambda\}.
\end{talign}
We cannot get a tight bound for~\eqref{eq:Seps_ub}, since this requires the (unknown) value function $V(\cdot)$.
Similar to the forward and backward steps, the saturation table allows us to remove the storage of $T$ separate saturated sets $S^k_t$ from EDDP and instead store a single look-up table, which reduces the memory cost by a factor of $T$.
}

\revone{
Having introduced $\Seps$, we can now observe that the auxiliary search point $z^k$ in Line~\ref{line:z_def} is the trial point from scenarios $i=1,\ldots,N$ with the largest upper bound on their saturation level.
Intuitively, the cutting-plane model may be the least accurate at $z^k$ among the trial points $\{x^k_i\}_{i=1}^k$.
We call this the auxiliary search point since it helps select the search point $x^k$ in Line~\ref{line:line8}.
}

In the next section, we examine the convergence properties of Inf-EDDP.

\subsection{Iteration complexity of Inf-EDDP}
First, we describe some relations between the cutting-plane model and the true value function, such as \revtwo{verifying that the cutting-plane model is an under-approximation of the true value function $V(\cdot)$}.
Although the following proof is similar to~\cite[Lemma 5]{lan2022complexity}, 
it mainly differs in that we do not know the true value function in the last stage since we are dealing with an infinite-horizon model, so the proof by induction goes forward in time rather than backward.
\begin{lemma} \label{lem:eddp_relation}
  For any iteration $k \geq 1$ and $i = 0,\ldots,N$,
 \begin{talign}
   \uV^{k-1}(x) \leq \uV^k(x) \leq \frac{1}{N} 
   \sum_{j=1}^N \tv_j^k(x) \leq V(x), &\ \forall x \in \mathcal{X}(\bar{\epsilon}) \label{eq:eddpV.1} \\
   \uF^{k-1}_i(x) \leq \uF^k_i(x) \leq F_i(x), &\ \forall x \in \mathcal{X}(\bar{\epsilon}), \label{eq:eddpF.1}
 \end{talign}
 where $F_0 := F$ is from~\eqref{eq:qgp_a3} and $\uF_0^k := \uF^k$ is from~\eqref{eq:policy}.
\end{lemma}
\begin{proof}
  Similar to~\cite[Lemma 5]{lan2022complexity},~\eqref{eq:eddpF.1} follows from~\eqref{eq:eddpV.1}, so we skip its proof here. 

  So we only need to show~\eqref{eq:eddpV.1}.
  The first relation $\uV^{k-1}(x) \leq \uV^k(x)$ is a direct consequence of the definition of $\uV^k$ in~\eqref{eq:uV_max}.
  
  To show the remaining inequalities, we proceed by induction on the iteration $k=0,\ldots$, where for notational convenience, we set \revtwo{$\uV^{-1}(x) = \uV^0(x) = \tilde v_j^0(x)$} for all $j =1,\ldots,N$. 
  To start, in view of the definition of $V$ in~\eqref{eq:qgp_a3b} and the lower bound $\underline h$ from Corollary~\ref{corr:h_minimizes}, then $\uV^0(x)$ as initialized in Line~\ref{line:line0} is a lower bound of $V(x)$. 
  So the base case goes through.
  We now consider the induction step.
  Let $x^{*}_i$ be the optimal solution for the optimization problem defined by \revtwo{$v_i(x)$ from~\eqref{eq:qgp_a3b}}, which is
  also feasible for the optimization problem defined by $\tv_i^k(x)$.
  Using the inductive hypothesis, 
  \begin{talign*}
    \frac{1}{N} \sum\limits_{i=1}^N \tv_i^k(x) 
    &\leq 
    \frac{1}{N} \sum\limits_{i=1}^N [h(x^{*}_i,\tilde c_i) + \lambda \uV^{k-1}(x^{*}_i)] 
    \leq \frac{1}{N} \sum\limits_{i=1}^N [h(x^{*}_i,\tilde c_i) + \lambda V(x^{*}_i)] = V(x).
  \end{talign*}
  
  We still need to show $\uV^k(x) \leq N^{-1}\sum_{i=1}^N \tv^k_i(x)$. 
  Observe that $\mathrm{cut}^k(x)$, which is defined as the second term in~\eqref{eq:uV_max}, is a supporting hyperplane centered at $x^{k-1}$.
  The validity of the supporting hyperplane comes from the fact that $V(\cdot)$ is convex by Lemma~\ref{lem:eddp_asmp_smooth1}. 
  Recalling that $\uV^k(x) = \max\{\underline{V}^{k-1}(x), \mathrm{cut}^k(x)\}$, we wish to show both terms in the aforementioned ``max'' lower bound $N^{-1}\sum_{i=1}^N \tv^k_i(x)$.
  Starting with $\uV^{k-1}$, we get   \begin{talign*}
    \uV^{k-1}(x) 
    \leq
    N^{-1} \sum_{i=1}^N \tv^{k-1}_i(x) 
    \leq 
    N^{-1} \sum_{i=1}^N \tv^{k}_i(x). 
  \end{talign*}
  where the first inequality is by the inductive hypothesis, and the second is by recalling the definition of $\tv_i^k$ and applying the inductive hypothesis for $\uV^{k-2}(u) \leq \uV^{k-1}(u)$ over all $u \in \mathcal X$.
  The supporting plane $\mathrm{cut}^k(x)$ lower bounds the convex value function $\tv_i^k$ (c.f.~Lemma~\ref{lem:eddp_asmp_smooth1}) via the gradient inequality.
  This completes the inductive hypothesis step and the proof too.
\end{proof}
Next, we will establish that $F_i$ and its lower bound $\uF^k_i$ are Lipschitz continuous.
Recall the parameter $\bar \epsilon \in (0,+\infty)$ satisfying Assumption~\ref{asmp:eddp_asmp_relint}, and the bounds $\overline{h},\underline{h}$ on the single-stage cost function $h(\cdot,\cdot)$ from Corollary~\ref{corr:h_minimizes}.
\begin{lemma} \label{lem:Fsmth.2}
  The function $F_i$ is Lipschitz continuous with constant $M \leq \revone{\frac{\overline{h} - \underline{h}}{(1-\lambda) \bar \epsilon}}$ for all $i=0,\ldots,N$, or
  \begin{talign*}
    |F_i(x) - F_i(x')| \leq M \|x-x'\|_\infty, \ \forall x,x' \in \mathcal X,
  \end{talign*}
  where $F_0 := F$ is from~\eqref{eq:qgp_a3}.
  Similarly, $\uF^k_i$ is Lipschitz continuous with constant $\underline M \leq \revone{\frac{\overline{h} - \underline{h}}{(1-\lambda) \bar \epsilon}}$ for every $k \geq 1$ and all $i=0,\ldots,N$, where $\uF_0^k := \uF^k$ is from~\eqref{eq:policy}.
\end{lemma}
\begin{proof}
  Proof of $F_i$ being Lipschitz continuous is directly from Corollary~\ref{corr:h_minimizes} and Lemma~\ref{lem:eddp_asmp_smooth1}.
  To show that $\uF^k_i$ is Lipschitz continuous for all $i \geq 0$, we need to show both its summands $h(\cdot, \tilde{\B c}_i)$ and $\uV^k(\cdot)$ are Lipschitz continuous.
  The former summand was shown in Corollary~\ref{corr:h_minimizes}. 
  To show the latter summand is Lipschitz continuous, we follow Lemma~\ref{lem:eddp_asmp_smooth1}, where we must show it is convex and proper.
  Convexity is immediate since the cutting-plane model $\uV^k$ is the pointwise max of linear functions, and hence it is a convex function~\cite[Theorem 5.5.2]{nesterov2018lectures}.
  To show it is proper, we will show it is absolutely bounded by $V(x) \in [\underline{h}/(1-\lambda), \overline{h}/(1-\lambda)]$ for any $x \in \mathcal{X}(\bar \epsilon)$.
  The lower bound is because of~\eqref{eq:eddpV.1} and the initialization $V^0(x) = \underline{h}/(1-\lambda)$, while the upper bound can be shown using the proof from Lemma~\ref{lem:eddp_asmp_smooth1} and $\uV^k(x) \leq V(x)$ from~\eqref{eq:eddpV.1}.
  Then we can apply the same argument from Lemma~\ref{lem:eddp_asmp_smooth1} to show $\uV^k(\cdot)$ is $M_{\uV}$-Lipschitz with $M_{\uV} \leq \frac{\overline{h}-\underline{h}}{(1-\lambda) \bar \epsilon}$.
  Thus, we have shown $\uF_i^k$ is Lipschitz continuous with constant $\underline{M} = M_h + \lambda M_{\uV}$.
\end{proof}

The next result provides sufficient conditions for a trial point to become saturated (Definition~\ref{def:dqgp_a0}). 
This result is similar to~\cite[Proposition 4]{lan2022complexity}, but our proof leverages the table $\mathcal{S}_\epsilon$.
Since $\mathcal{S}_\epsilon$ changes between iterations, we say $\mathcal{S}_\epsilon^k$ is the table $\mathcal{S}_\epsilon$ at the beginning of iteration $k$, or equivalently, at the end of iteration $k-1$. 
\begin{proposition} \label{prop:eddp_next_sat}
  Let $\epsilon_t$ follow the recursive definition $\epsilon_{t-1} = (M + \underline M)\epsilon + \lambda \epsilon_t$ and $\epsilon_{T}$ be defined so that for any iteration, every point in $\mathcal{X}$ is $\epsilon_{T}$-saturated at any iteration $k > 0$. 
  Suppose $\epsilon \in (0,+\infty)$ is selected s.t.~$\epsilon_t$ is nondecreasing with $t$.
  Then for any $k \geq 1$, the following two conditions hold:
  \begin{enumerate}[(a)]
      \item \textit{Case (a)} For any $x \in \mX$, if $\Seps^k (\revone{\Geps}(\revtwo{x})) = t$, then for all $i=0,\ldots,N$,
           \begin{talign} \label{eq:F_bound_eps}
             \lambda [ V(\revtwo{x}) - \uV^{k-1} (\revtwo{x}) ] 
             =
             F_i(\revtwo{x}) - \uF_i^k (\revtwo{x}) 
             \leq
             \epsilon_{t-1},
           \end{talign}
           where $F_0 := F$ is from~\eqref{eq:qgp_a3} and $\uF_0^k := \uF^k$ is from~\eqref{eq:policy}.
           \sloppy \revtwo{Furthermore, when $\max_{i =1,\ldots,N} \Seps^k (\revone{\Geps}(x^k_i)) = t$, where the trial point $x^k_i$ is from Line~\ref{line:line2}, then $x^{k-1}$ is $\epsilon_{t-1}$-saturated at iteration $k$, or 
           \begin{talign} \label{eq:next_stage_sat}
              V(x^{k-1}) - \uV^k(x^{k-1}) \leq \epsilon_{t-1}.
           \end{talign}}
        \item \textit{Case (b)} For any $x \in \mathcal X$, if $\Seps^{k+1}(\revone{\Geps}(x)) = t < T$, then there exists a previous trial point $x^j$ (where $j \leq k$) s.t.~$\|x-x^j\|_\infty \leq \epsilon$ and $x^j$ is $\epsilon_{t'}$-saturated for any $t' \in [t,T]$ at iteration $k+1$.
  \end{enumerate}
\end{proposition}
\begin{proof}
    First, we note that the equality in~\eqref{eq:F_bound_eps} holds by definition of $\uF_i$. 

    Now, we prove the remainder of Case (a) and Case (b) by induction on $k \geq 0$. 
    For the base case $k=0$, Case (a) is vacuously true since Inf-EDDP starts at iteration $k=1$.
    Similarly, Case (b) is vacuously true since we initialize $\Seps(\revone{\Geps}(x)) = T$ (Line~\ref{line:line1}) before iteration $1$.

    We now consider the inductive case.
    Starting with Case (a), by invoking the inductive hypothesis of Case $(b)$ towards the standing assumption $\Seps^k(\revone{\Geps}(\revtwo{x})) = t$, 
    then we ensure there is a previous search point $x^{j}$ (where $j < k$) that is $\epsilon_{t}$-saturated and $\|x^{j} - x\|_\infty \leq \epsilon$.
    Using an argument similar to~\cite[Proposition 4]{lan2022complexity}, we have $F_i(x) - \uF_i^k(x)
        \leq
        (M + \underline{M})\epsilon + \lambda \epsilon_{t}
        =
        \epsilon_{t-1}$.
    This establishes~\eqref{eq:F_bound_eps} for the inductive case.
    Now, to show~\eqref{eq:next_stage_sat}, we first recall $x^k_i$ is an optimal solution to $\tv^k_i(x^{k-1})$ from~\eqref{eq:eq5} for each $i=1,\ldots,N$.
    Then $x^k_i$ is a feasible solution to the problem described by $\revone{v_i}(x^{k-1})$, hence
    \begin{talign} \label{eq:qgp_a43}
        V(x^{k-1})
        =
        N^{-1}\sum_{i=1}^N v_i(x^{k-1})
        \leq
        N^{-1}\sum_{i=1}^N F_i(x^k_i).
    \end{talign}
    Furthermore, using the definition of $\uV^k$,
    \begin{talign} \label{eq:qgp_a42}
        \uV^k(x^{k-1})
        \geq
        N^{-1}\sum_{i=1}^N \tv^k_i(x^{k-1})
        =
        N^{-1}\sum_{i=1}^N \uF^k_i(x^k_i).
    \end{talign}
    Combining the above two with~\eqref{eq:F_bound_eps} yields~\eqref{eq:next_stage_sat}.

    Now to finish, we consider the inductive case for Case (b) for iteration $k$.
    Fix an arbitrary $x \in \mathcal X$.
    In view of Line~\ref{line:line5}, $\Seps(\revone{\Geps}(x))$ strictly decreases to $t-1$ during iteration $k$ if and only if both
    \begin{enumerate}[1),nosep]
        \item $\Seps(\revone{\Geps}(x^{k-1})) > t := \Seps^k(\revone{\Geps}(z^k))$
        \item $\revone{\Geps(x) = \Geps(x^{k-1})}$.
    \end{enumerate}
    If one of these \revtwo{does} not hold, then we instead have $\Seps^{k+1}(\revone{\Geps}(x)) = \Seps^k(\revone{\Geps}(x)),~\forall x \in \mX$.
    Invoking the inductive hypothesis, this shows Case $(b)$ holds for $x$ at iteration $k$.
    Otherwise, if both 1) and 2) occur, we claim $x^{k-1}$ is the previous trial point that satisfies Case $(b)$ for $x$.
    To prove the claim, we use the fact that $\revone{\Geps(x) = \Geps(x^{k-1})}$ to ensure there is a $u \in \mXeps$ such that $\|x-x^{k-1}\|_\infty \leq \|x-u\|_\infty + \|u-x^{k-1}\|_\infty \leq \frac{\epsilon}{2} + \frac{\epsilon}{2} = \epsilon$, where the first inequality is by construction of the $\epsilon$-net. 
    Since $t = \Seps^k(\revone{\Geps}(z^k)) = \max_{i=1,\ldots,N} \Seps^k(\revone{\Geps}(x^k_i))$, then we can apply~\eqref{eq:next_stage_sat} for iteration $k$ to verify that $x^{k-1}$ is $\epsilon_{t-1}$-saturated at iteration $k$, and hence also $\epsilon_{t'}$-saturated at iteration $k$ for all $t' \in [t-1,T]$ by the nondecreasing property of $\epsilon_t$.
\end{proof}

Unlike the finite-horizon case, the infinite-horizon \revone{case} does not have access to an exact last-stage value function because there is no last stage in the infinite-horizon case.
Instead, we will show every feasible point is $\epsilon_{T}$-saturated for some $\epsilon_{T} > 0$ that depends on the bounds $\overline{h}, \underline{h}$ of the single-stage cost function $h(\cdot,\cdot)$ from Corollary~\ref{corr:h_minimizes}. 
\begin{lemma} \label{lem:eddp_every_pt_is_sat}
  For any $k \geq 0$ and $x \in \mathcal{X}$, then $V(x) - \uV^k(x)  \leq  (1-\lambda)^{-1}[\overline{h}-\underline{h}]$.
\end{lemma}
\begin{proof}
  Let $\hat{x}_i \in \mathrm{argmin}_{x \in \mathcal X(\bar \epsilon)} h(x, \tilde{\B c}_i)$ 
  and $\bar{x}_i$ be the solution to the optimization problem $v_i(x)$ as defined in~\eqref{eq:qgp_a3}.
  We have
  \begin{talign*}
    V(x) - \uV^k(x)
    &\leq 
    V(x) - \uV^0(x) \\
    &= 
    N^{-1} \sum_{i=1}^N \big( [h(\bar{x}_i, \tilde{\B c}_i) - h(\hat{x}_i, \tilde{\B c}_i)] + \lambda[V(\bar{x}_i) - \uV^0(\hat{x}_i)] \big) \\
    &\leq
    N^{-1} \sum_{i=1}^N \big( [\overline{h}-\underline{h}] + \lambda[V(\bar{x}_i) - \uV^0(\hat{x}_i)] \big) \\
    &=
    [\overline{h}-\underline{h}] + {\lambda} {N}^{-1} \sum_{i=1}^N [V(\bar{x}_i) - \uV^0(\bar{x}_i)],
  \end{talign*}
  where the first line is by Lemma~\ref{lem:eddp_relation}, 
  the second line is by definition of $V$ and $\uV^0$,
  the third line is Corollary~\ref{corr:h_minimizes},
  and finally the last line is by the fact that $\uV^0$ is a constant function.
  Recursively applying the above result and using the fact that $V$ and $\uV^0$ are both bounded (the former can be shown by Lemma~\ref{lem:eddp_asmp_smooth1} while the latter by Line~\ref{line:line0} and finiteness of $\underline h$) yields the result.
\end{proof}

We next make an important observation regarding how the table $\mathcal{S}_\epsilon$ changes during the algorithm.
This lemma is adapted from~\cite[Proposition 5]{lan2022complexity} to make use of the table $\Seps$.
And recall that $\mathcal{S}_\epsilon^k$ is the table $\mathcal{S}_\epsilon$ at the beginning of iteration $k$, or equivalently, the end of iteration $k-1$. 
\begin{lemma} \label{lem:eddp_find_or_converge_first}
  Let $k = \ell T+1$ for an integer $\ell \geq 0$. Inf-EDDP either terminates
  during iteration $k$ or generates a search point $\revone{x}^{k+t-2}$ 
  satisfying $\Seps^{k+t-2}(\revone{\Geps}(\revone{x}^{k+t-2})) \geq t$ and $\Seps^{k+t}(\revone{\Geps}(\revone{x}^{k+t-2})) \leq t-1$ for some $t$ where $t \in [2,T]$. 
\end{lemma}
\begin{proof}
    First, we observe that $x^{k+t-1} = z^{k+t-1}$ holds for all $t \in [2,T]$ by the update in Line~\ref{line:line8}.
    Second, we observe that $\Seps(\Geps(x)) \leq T,~\forall x \in \mX$ holds in every iteration of Inf-EDDP, because we initialize $\Seps(\cdot) = T$ in Line~\ref{line:line0}, and all updates to $\Seps$ in Line~\ref{line:line5} are non-increasing.
    In particular, we have $\Seps(\Geps(z^{k+T-1})) \leq T$.
    With these two observations, we can deduce that one of the following $T$ cases must occur:
    \begin{enumerate}
        \item Case 1: $\Seps^k(\revone{\Geps}(x^k)) \leq 1$
        \item Case $t \in [2, T]$: $\Seps^{k+t-2}(\revone{\Geps}(\revone{x}^{k+t-2})) \geq t$ and $\Seps^{k+t-1}(\revone{\Geps}(z^{k+t-1})) \leq i, 
        \ \forall i \in [t,T]$.
    \end{enumerate}
    We start with Case 1. 
    Note that $x^k = x^k_0$ because of $k = \ell T+1$ and Line~\ref{line:line8}.
    As a result, the assumption of Case 1 implies $\Seps^k(\revone{\Geps}(x^k_0)) \leq 1$, which ensures Inf-EDDP terminates during iteration $k$ at Line~\ref{line:line6}.

    We now consider Case $t$. 
    \sloppy The second half of the hypothesis for Case $t$ says $\Seps^{k+t-1}(\revone{\Geps}(z^{k+t-1})) \leq t$ for some $t \in [2,T]$.
    Using this observation with the update in Line~\ref{line:line5} ensures that after iteration $k+t-1$ the saturation table satisfies $\Seps^{k+t}(\revone{\Geps}(x^{k+t-2})) \leq t-1$.
    Combining this with the first half of the hypothesis for Case $t$ completes the proof.
\end{proof}
We are now ready to establish the number of iterations of Inf-EDDP to output a solution $x^k_0$ that approximately solves~\eqref{eq:qgp_a3}.
\revtwo{
We define the size of $\Seps$ as
\begin{talign} \label{eq:S_size}
    \vert \Seps \vert
    :=
    \sum_{x \in \mathcal \mXeps} \max\{0, \Seps(x)\}.
\end{talign}
}
Let $F^* = \min_{x \in X(x^0,\xi^1)} F(x)$ be the optimal value of the problem of interest~\eqref{eq:qgp_a3}. 

\begin{theorem} \label{thm:thm2}
  For any $T \geq 1$ and $\epsilon \in \big(0, \frac{\lambda^T(\overline{h}-\underline{h})}{M + \underline{M}} \big]$, let $\epsilon_t$ be \revone{recursively defined as in Proposition~\ref{prop:eddp_next_sat}}.
  Then, Inf-EDDP returns a feasible solution $x^k_0$ such that
  \begin{talign} \label{eq:qgp_a7}
    F(x^k_0) - F^* \leq \epsilon_0 
    \leq 
    \frac{(M+\underline{M})\epsilon}{1-\lambda} + \frac{\lambda^{T}(\overline{h}-\underline{h})}{1-\lambda}
  \end{talign}
  in at most $K$ iterations, where $K := T^2\left(D/\epsilon+1\right)^{n}$.
\end{theorem}

\begin{proof}
  We start by bounding the number of \revtwo{iterations}.
  We recall that the $\epsilon$-net $\mathcal X_\epsilon$ of a bounded set $\mathcal X$ with diameter $D$ and dimension $n$ has size at most $\vert \mathcal X_\epsilon \vert \leq \bar{K} := ({D}/{\epsilon}+1)^n$.
  Therefore, from the definition of the size of $\mathcal{S}_\epsilon$ from~\eqref{eq:S_size} and the fact that we initialize all points $x \in \mathcal X$ with $\mathcal{S}_\epsilon(x) = T$, then $\vert \mathcal{S}_\epsilon^1 \vert \leq T\bar{K}$.
  Moreover, by definition of the size of $\Seps$ in~\eqref{eq:S_size}, the condition $\vert \mathcal{S}_\epsilon^k \vert \leq 1$ implies $\Seps(x) \leq 1, \ \forall x \in \mathcal X$. 
  In particular, by fixing $x=x^k_0$, we find $\Seps(x^k_0) \leq 1$, which means Inf-EDDP terminates by \revtwo{Line~\ref{line:line6}} during iteration $k$.
  Now, we claim that during epoch $\ell \geq 0$, which is defined as iterations $[\ell T+1, (\ell+1)T]$, either (1) Inf-EDDP terminates or (2) $\vert \mathcal{S}_\epsilon \vert$ decreases by at least one.
  If this claim is true, then we observe that (2) can only hold for at most $T\bar{K}$ epochs (i.e., $T^2 \bar{K}$ iterations) since after epoch $T\bar{K}$, we would have $\vert \mathcal{S}^{T^2 \bar{K}}_\epsilon \vert \leq 1$.

  To prove the claim, we first define $k = \ell T + 1$ and $t' = k+t-2$ for some $t \in [2,T]$.
  Hence, $t' \in [\ell T+1, (\ell+1)T-1]$.
  We use Lemma~\ref{lem:eddp_find_or_converge_first} to show that either Inf-EDDP terminates at iteration $t'$ or there is a search point $x^{t'}$ satisfying $\Seps^{t'}(\revone{\Geps}(x^{t'})) \geq t$ and $\Seps^{t'+2}(\revone{\Geps}(x^{t'})) \leq t-1$. 
  Because $\revone{\Geps}(x^{t'}) \in \mXeps$ and the update to $\mathcal{S}$ are monotone decreasing (Line~\ref{line:line5}), the size $\vert \Seps \vert$ decreases by at least $\Seps^{t'}(\revone{\Geps}(x^{t'})) - \Seps^{t'+2}(\revone{\Geps}(x^{t'})) \geq 1$ during iterations $t'$ and $t'+1$, as desired.

  \revone{The bound on $\epsilon_0$ can be deduced from the recursive definition, $\epsilon_{t-1}=(M + \underline{M})\epsilon + \lambda \epsilon_t$ for $t \leq T$ and $\epsilon_T = [\overline{h}-\underline{h}]/(1-\lambda)$ from Proposition~\ref{prop:eddp_next_sat}.}
  The bound~\eqref{eq:qgp_a7} is due to the fact when the algorithm terminates (Line~\ref{line:line6}), then $\mathcal{S}_\epsilon(x^k_0) \leq 1$, which ensures
  \begin{talign} \label{eq:F_optimality_relationship}
      F(x^k_0) - F(x^*)
      \leq 
      F(x^k_0) - \uF^k_0(x^*)
      \leq
      F(x^k_0) - \uF^k_0(x^k_0)
      \leq
      \epsilon_0,
  \end{talign}
  where the first inequality is by Lemma~\ref{lem:eddp_relation}, the second since $x^k_0$ is the minimizer of $\uF^k_0(\cdot)$, and the last inequality by~\eqref{eq:F_bound_eps}.
  By choice of $\epsilon$ and definition of $\epsilon_t$, one can verify $\epsilon \leq \frac{(1-\lambda)\epsilon_t}{M + \underline{M}}$ for all $t$. 
  In view of the recursive definition of $\epsilon_{t-1} = (M + \underline{M})\epsilon + \lambda\epsilon_t$, we can ensure that $\epsilon_t$ is nondecreasing with $t$.
  Then we can apply Proposition~\ref{prop:eddp_next_sat} to ensure~\eqref{eq:qgp_a7}.
\end{proof}
A few remarks are in order. 
The theorem shows a trade-off between accuracy and complexity when choosing the effective planning horizon $T$.
In particular, the accuracy bound $F(x^k_0) - F^* \leq \varepsilon$ can be achieved by setting $T =\ln(\frac{2(\overline{h}-\underline{h})}{\varepsilon(1-\lambda)})/(1-\lambda)$ and $\epsilon = \frac{\lambda^T(\overline{h}-\underline{h})}{2(M + \underline{M})}$.
Now, for a quick comparison, we observe that applying EDDP to a truncated $T$-horizon problem~\cite{lan2022complexity} requires $O(T^2[D/\epsilon+1]^n)$ forward and backward steps, where each forward and backward phase has $T$ forward and $T-1$ backward steps, respectively.
This matches the iteration complexity of Inf-EDDP, where every iteration requires only one forward and one backward step.

In the next section, we will present a new search point selection rule that improves the iteration complexity of Inf-EDDP by a factor of $T$, and thus provides a clear advantage over directly applying EDDP.

\subsection{Reducing the dependence on the effective planning horizon} \label{sec:better_T}
\revone{
  We propose a new search point selection that allows Inf-EDDP to explore the feasible region longer.
  As seen in Algorithm~\ref{alg:alg5}, the auxiliary search point $z^k$ is replaced by an augmented auxiliary search point $z^k_+$, which additionally considers the here-and-now solution $x^k_0$ as a potential search point with the largest saturation level.
  The algorithm also explores the feasible region longer with a trajectory length doubled to $2T$.
  These changes allow for more adaptive exploration by resetting the search point $x^k$ to $x^k_0$ when the latter is under-explored. 
  Otherwise, the algorithm continually explores for a longer period of $2T$.
  Therefore, we call the algorithm continually-exploring Inf-EDDP, or CE-Inf-EDDP.}

\begin{algorithm}[h!]
\caption{}
\label{alg:alg5}
\begin{algorithmic}[1]
\Procedure{Continually Exploring Inf-EDDP}{$\epsilon, T$}
    \State{Run Inf-EDDP (Algorithm~\ref{alg:alg2}) but replace $z^k$ in Lines~\ref{line:z_def} and~\ref{line:line5} with
      \begin{talign*} 
        z^k_+ := \mathrm{argmax}_{\{x^k_i\}_{i=0}^N} \Seps(\Geps(x_i^k)),
      \end{talign*} 
      and replace Line~\ref{line:line8} with
      \begin{talign*}
        x^k 
        =
        \begin{cases}
          x^k_0 ~&: ~ k~\mathrm{mod}~(2T) = 1 \\
          z^k_+ ~&: ~ \text{o.w.}
        \end{cases}.
      \end{talign*}
    } \label{line:next_search_modify}
\EndProcedure
\end{algorithmic}
\end{algorithm}

\revone{Note that doubling the trajectory length to $2T$ is arbitrary, and our theoretical results are nearly identical when the trajectory length is increased to $\lceil c \cdot T \rceil$ for a $c \geq 2$. 
Empirically, $c=2$ performs better since a larger $c$ can cause the algorithm to over-explore and rarely explore the feasible region at the here-and-now solution $x^k_0$.
}

We now analyze the convergence of the CE-Inf-EDDP algorithm. 
\revtwo{One can check that Lemmas~\ref{lem:eddp_relation},~\ref{lem:Fsmth.2}, and~\ref{lem:eddp_every_pt_is_sat} still hold, since their proofs do not depend on the selection of $z^k$. 
Proposition~\ref{prop:eddp_next_sat} is also still true, and the proof is nearly identical with the augmented auxiliary point $z^k_+$. 
In contrast, the results for Lemma~\ref{lem:eddp_find_or_converge_first} and Theorem~\ref{thm:thm2} differ based on the selection of $z^k$.
First, we adapt Lemma~\ref{lem:eddp_find_or_converge_first} by incorporating this longer exploration strategy.
}
\begin{lemma} \label{lem:eddp_find_or_converge}
  Consider a $k \in [2\ell T+1, (2\ell+1)T+2]$ for any integer $\ell \geq 0$.
  CE-Inf-EDDP either terminates during iteration $k$ or generates a search point $x^{k+t-2}$ satisfying $\mathcal{S}_\epsilon^{k+t-2}(x^{k+t-2}) \geq t$ and $\mathcal{S}_\epsilon^{k+t}(x^{k+t-2}) \leq t-1$ for some $t$ where $t \in [2,T]$. 
\end{lemma}
\begin{proof}
  Using a similar proof to Lemma~\ref{lem:eddp_find_or_converge_first}, we can deduce that one of the following $T$ cases must occur:
    \begin{enumerate}
        \item Case 1: $\Seps^k(\revone{\Geps}(x^k)) \leq 1$
        \item Case $t \in [2, T]$: $\Seps^{k+t-2}(\revone{\Geps}(x^{k+t-2})) \geq t$ and $\Seps^{k+t-1}(\revone{\Geps}(x^{k+t-1})) \leq i, 
        \ \forall i \in [t,T]$.
    \end{enumerate}
    We start with Case 1. 
    If $k=2\ell T +1$, then we get $x^k = x^k_0$ by the modified Line~\ref{line:next_search_modify}. 
    This implies $\Seps^k(\revone{\Geps}(x^k_0)) \leq 1$ based on the hypothesis Case 1, and then the algorithm terminates during iteration $k$ at Line~\ref{line:line6}.  
    Otherwise, if $k \ne 2 \ell T+1$ for some $\ell$, then one can similarly show that $\Seps^k(\revone{\Geps}(x^k)) = \Seps^k(\revone{\Geps}(z^k_+)) = \max_{i =0,\ldots,N} \Seps^k(\revone{\Geps}(x^k_i))$, which implies $\Seps^k(\revone{\Geps}(x^k_0)) \leq 1$. 
    Thus, CE-Inf-EDDP terminates during iteration $k$.

    We now consider Case $t \in [2,T]$. 
    Because $x^{k+t-1} = z^{k+t-1}_+$ due to the modified Line~\ref{line:next_search_modify}, the second half of the hypothesis for Case $t$ ensures that $\Seps^{k+t-1}(\revone{\Geps}(z^{k+t-1}_+)) \leq t$.
    Using this observation with the update in the modified Line~\ref{line:line5} (where $z^k$ replaces $z^k_+$), we find that after iteration $k+t-1$, the saturation table satisfies $\Seps^{k+t}(\revone{\Geps}(x^{k+t-2})) \leq t-1$.
    Combining this inequality with the first half of the hypothesis for Case $t$ completes the proof.
\end{proof}

Now, the above lemma can be used as follows, which is the main insight into improving the dependence $T$. 
Invoking the above lemma on each of the $T+1$ iterations in $[2\ell T+1, (2\ell+1)T+2]$ ensures the size $\vert \Seps \vert$ decreases by $O(T)$ during this period.
In contrast, directly applying Lemma~\ref{lem:eddp_find_or_converge_first} can only guarantee a decrease of at least $2$ within this same period.
It seems more exploration can improve the accuracy of the cutting-plane model.
We now formalize this below, where the proof only looks at odd or even indices to avoid double counting. 
\begin{lemma} \label{lem:eddp.1}
  Consider iterations $[2\ell T+1, 2(\ell+1) T]$ for some $\ell \geq 0$.
  Then CE-Inf-EDDP decreases $\vert \Seps \vert$ by at least $\lceil (T+2)/2 \rceil$ or terminates during this period.
\end{lemma}
\begin{proof}
  For any iteration $k \in I_\ell := [2\ell T+1, (2\ell +1)T+2]$, let $t$ be the index where iteration $\tau(k) := k+t-2$ satisfies the properties of Lemma~\ref{lem:eddp_find_or_converge}. 
  We call $\tau(k)$ the \textit{target iteration} of $k$. 
  Split $I_\ell = I_e \dot{\cup} I_o$, where $I_e$ is the set of iterations whose target iteration is even, and $I_o$ is similarly defined but for odd-indexed target iterations.
  Let $I$ be the larger one of these two sets, and let $J := \{\tau(k)\}_{k \in I}$ be the corresponding set of target iterations.
  Since $\vert I_\ell \vert = T+2$, it must be that $\vert I \vert \geq (T+2)/2$.

  Note that multiple iterations in $I$ can have matching target iterations.
  With that in mind, for any target iteration $j \in J$, there exists a set of $L$ iterations $\{k_i\}_{i=1}^L \subseteq I$ that have the same target iteration, $j=k_i+t_i-2$, where $L$ is some positive integer satisfying $1 \leq L \leq T-1$, and $\{t_i\}_{j=1}^L  \subseteq \{2,\ldots,T-1\}$ is a set of integers that ensure $k_i+t_i-2=j$. 
  Without loss of generality, suppose $\{t_i\}_{i=1}^L$ is decreasing so $t_1$ and $t_L$ are the largest and smallest values from the set, respectively. 
  Invoking Lemma~\ref{lem:eddp_find_or_converge} to iterations $\{k_i\}_{i=1}^L$ (and recalling the assumption about not terminating), we find $\Seps^{\revone{k_1+t_1-2}}(\revone{\Geps(x^{k_1+t_1-2})}) \geq t_1 $ and $ \Seps^{\revone{k_L+t_L}}(\revone{\Geps(x^{k_L+t_L-2})}) \leq t_L-1$. 
  Equivalently, $\Seps(\revone{\Geps}(\revone{x^j}))$ decreased by at least $t_1-t_L+1$ for a search point $x^j$ during target iteration $j$ and iteration $j+1$. 
  Since $L \leq t_1 - t_L + 1$, because the set $\{t_j\}_{j=1}^L$ contains strictly decreasing integers, we can conclude that $\Seps(\revone{\Geps}(x^j))$ decreased by at least $L$ during those two iterations. 
  Repeating this argument for every $j \in J$ ensures that $|\Seps|$ decreased by $\vert I \vert$ during the mentioned period, and we recall that $\vert I \vert > (T+2)/2$.
\end{proof}

The above lemma asserts that we can, on average, decrease $\vert \mathcal{S}_\epsilon\vert$ by at least $1$ every $4$ iterations rather than $1$ every $T$ iterations. 
This more frequent decrease of $\vert \Seps \vert$ leads to improved iteration complexity.
\begin{theorem} \label{thm:thm2b}
  Let everything be defined as in Theorem~\ref{thm:thm2}. CE-Inf-EDDP returns a feasible solution $x^k_0$ such that
  \begin{talign*} 
    F(x^k_0) - F^* 
    \leq
    \epsilon_0
    \leq
    \frac{(M+\underline{M})\epsilon}{1-\lambda} + \frac{\lambda^{T}(\overline{h}-\underline{h})}{1-\lambda},
  \end{talign*}
  but now the number of iterations is at most $K := 4T (D/\epsilon+1)^{n}$.
\end{theorem}
\begin{proof}
  The proof is largely the same as Theorem~\ref{thm:thm2}, except we apply Lemma~\ref{lem:eddp.1} to show $\vert \Seps \vert$ must either decrease by at least $\lceil T/2 \rceil$ every $2T$ iterations.
\end{proof}
Compared to Theorem~\ref{thm:thm2}, the number of iterations is reduced by $T/4$. 
While a similar improvement was shown for the more general finite-horizon case~\cite{zhang2020distributionally}, the modifications are more complicated as they require both a lower- and upper-bound model of the value function.
In contrast, CE-Inf-EDDP only requires a lower-bound model, which simplifies the computational cost of each iteration. 
In the next section, we show that CE-Inf-EDDP can handle alternative search point selections.

\section{Other trial point selection strategies} \label{sec:sec4}
Our goal in this section is to extend some of the trial point selection rules (Line~\ref{line:next_search_modify}) from the existing literature, such as using a duality gap~\cite{baucke2017deterministic,zhang2020distributionally,zhang2022stochastic} or random sampling~\cite{pereira1991multi,homem2011sampling,shapiro2011analysis}, towards the infinite-horizon setting.

\subsection{Selection by Largest Gap via Dual Bound} \label{sec:upperbounds}
Recall the cutting-plane model $\uV$ provides a lower bound of the value function $V$ (Lemma~\ref{lem:eddp_relation}).
We can also develop \revtwo{an} upper-bound model $\oV$.
Once constructed, it can be used for two purposes: first, to derive a duality gap and possibly terminate the algorithm sooner, and second, to use the aforementioned duality gap as a trial point selection \revone{criterion}.

To start, fix some $\overline{M}_0 \in [0,+\infty)$. 
We build the upper bound recursively by
\begin{talign} \label{eq:qgp_a15}
  \oV^k(x) = 
  \begin{cases}
    (1-\lambda)^{-1} \overline{h} &: k = 0 \\
    {N}^{-1} \sum_{i=1}^N \ov^k_i(x) &: k \geq 1
  \end{cases},
\end{talign}
where for some set of previous trial points $\{x^j \in \mathcal X\}_{j=1}^k$,
the $i$-th scenario's model is
\begin{align} \label{eq:upper_V} 
  \ov^k_i(\bar x) :=  \max_{\mu,\rho} \ & \mu + \langle \rho, \bar x \rangle \\
  \text{s.t. } 
  &\mu + \langle \rho, x^{j} \rangle \leq \hat{v}^j_i(x^{j}), \ \forall j \in \{1,\ldots,k\} \nonumber \\
  & \|\rho\|_1 \leq \overline{M}_0, \nonumber
\end{align}
and
\begin{talign} \label{eq:qgp_a14}
  \hv^{k}_i(\bar x) := \min_{x \in X(\bar x, \tilde{\xi}_i)} \big \{ \oF^k_i(x) := h(x, \tilde{c}_i) + \lambda \oV^{k-1}(x) \big\}.
\end{talign}

Let us make some remarks.
First, the upper-bound model is initialized at iteration $k=0$, similarly to the lower bound (Line~\ref{line:line0}), where the cost upper bound is from Corollary~\ref{corr:h_minimizes}. 
Second, the two functions $\oF^k_i$ and $\oV^k$ are the upper-bound counterparts of $F_i$ and $V$, respectively, from~\eqref{eq:qgp_a3}.
Third, one can check that the maximization problem in~\eqref{eq:upper_V} is always feasible and bounded (since $\rho$ is bounded).
Fourth, the optimization problem~\eqref{eq:upper_V} is derived by taking the dual of the linear program (LP) that upper bounds scenario $i$'s value function $v_i(\bar{x})$ by taking the minimum value of the convex hull of previous trial points' values. 
The bound on the norm of $\rho$ is added to the dual LP to ensure the solution is finite.
See~\cite{baucke2017deterministic} for more details.

We now verify that the proposed functions are indeed upper bounds. 
Recall $\overline{M}_0$ from~\eqref{eq:upper_V} is a tunable hyperparameter.
\begin{lemma} \label{lem:lqgp_a2}
  The function $\ov^k_i(\cdot)$ is convex and $\overline{M}_0$-Lipschitz w.r.t.~the
  $\|\cdot\|_\infty$ norm. 
  In addition, if $v_i(x)$ from~\eqref{eq:qgp_a3} is $M_v$-Lipschitz over $\mathcal{X}$ for all $i$ (w.r.t.~the $\ell_\infty$-norm) and $\overline{M}_0 \geq M_v$, then
  $F_i(x) \leq \oF^k_i(x) \ \forall k \geq 1, x \in \mathcal X$.
\end{lemma}
\begin{proof}
  Convexity is immediate since $\bar{v}^k_i(\cdot)$ is the pointwise max of convex functions,
  while Lipschitz continuity is by~\cite[Lemma 2.1]{baucke2017deterministic}.

  To prove \revtwo{$\oF^k_i(\cdot)$} is an upper bound of $F_i(\cdot)$, it suffices to prove for all $k\geq 0$ and $x \in \mathcal X$,
  $v_i(x) \leq \ov^k_i(x)$, which implies $V(x) \leq \oV^k(x)$.
  Towards that, define an exact upper-bound model ${\ov^*}^k_i(\cdot)$,
  which is equivalent to $\ov^k_i(\cdot)$ except that the constraints' left-hand side values $\hv^j_i(x^j)$ from~\eqref{eq:upper_V} are replaced by the true value $v_i(x^j)$ from~\eqref{eq:qgp_a3}.
  Since $v_i(\cdot)$ is assumed to be $M_v$-Lipschitz and to satisfy $\overline{M}_0 \geq M_v$, then ${\ov^*}^k_i(x) \geq v_i(x)$ for any $x \in \mathcal{X}$~\cite[Lemma 2.4]{baucke2017deterministic}.
  Thus, all that we need to show is $\ov^k_i(x) \geq {\ov^*}^k_i(x)$.
  In view of~\eqref{eq:upper_V}, it suffices to show $\hv_i^j(x^j) \geq v_i(x^j)$ for all $j=1,\ldots,k$.
  This can be shown by a simple mathematical induction, so we skip this step to keep the proof concise.
\end{proof}
We make some remarks.
Note that the aforementioned assumption that $v_i(\cdot)$ is Lipschitz can be satisfied; see Lemma~\ref{lem:eddp_asmp_smooth1} for more details.
Additionally, one can establish a relationship between the upper bound models and true functions, similar to Lemma~\ref{lem:eddp_relation}.
The result is straightforward upon inspecting~\eqref{eq:upper_V} and using the above lemma, so we do not explicitly state the result.
Finally, in view of the definition of $\oF^k_i$ and the above lemma, we can ensure $\oF^k_i$ is $\overline{M}$-Lipschitz (w.r.t.~$\ell_\infty$-norm) over $\mathcal{X}$ with constant $\overline{M} := \lambda \overline{M}_0 + M_h$.

Now, equipped with both an upper and lower bound, we define the gap function,
\begin{align} \label{eq:qgp_a13}
  \gamma^{k-1}(x) = \oV^{k-1}(x) - \uV^{k-1}(x).
\end{align}
Notice that if $\gamma^{k-1}(x) \leq \epsilon$ for some $\epsilon \geq 0$ and $x \in \mathcal X$, then $V(x) - \uV^k(x) \leq \epsilon$.
So one can use the (computable) gap function to bound the (uncomputable) gap $V(x) - \uV^k(x)$.
With the gap function in place, we modify Inf-EDDP to select points using $\gamma^{k-1}$, resulting in Algorithm~\ref{alg:alg4}. 
The main difference is that we update $\mathcal{S}_\epsilon$ immediately after we solve each subproblem.  
When assessing the accuracy levels $\{\epsilon_t\}$ (given in Proposition~\ref{prop:pqgp_a1}), we assume they are strictly increasing, i.e.,  $\epsilon_{t-1} < \epsilon_t$. 
Now, for this section only, we say a point $x^{k}$ is $\epsilon_t$-saturated when 
\begin{talign*}
  \gamma^{k-1}(x^k) \leq \epsilon_t,
\end{talign*}
c.f. Definition~\ref{def:dqgp_a0}. 
Notice in the algorithm, we set $\epsilon_{-1} = -1$. This is an arbitrary negative number, which is for notational convenience when a trial point becomes $\epsilon_0$-saturated.

\begin{algorithm}[h]
\begin{algorithmic}[1]
\Procedure{Gap-Inf-EDDP}{$\epsilon, T$}
    \State{Modify Inf-EDDP (Algorithm~\ref{alg:alg2}) like in Algorithm~\ref{alg:alg5} and also add after Line~\ref{line:line2}
      \begin{align} \label{eq:qgp_a46}
          \mathcal{S}_\epsilon(x^k_i) = \min\big( \mathcal{S}(x^k_i), \{t : \epsilon_{t-1} < \gamma^{k-1}(x^k_i) \leq \epsilon_t\} \big) \ \text{ (where $\epsilon_{-1} = -1$)}.
      \end{align}
    } 
\EndProcedure
\end{algorithmic}
\caption{}
\label{alg:alg4}
\end{algorithm}

The following result is nearly identical to Proposition~\ref{prop:eddp_next_sat}. 
The main difference is that we replace the value functions with their upper bounds, and we account for the update to $\mathcal{S}_\epsilon$ in Gap-Inf-EDDP.
Hence, we skip the proof.

\begin{proposition} \label{prop:pqgp_a1}
  Let $\epsilon_t$ follow the recursive definition $\epsilon_{t-1} = (\overline M + \underline M)\epsilon + \lambda \epsilon_t$ and $\epsilon_{T}$ be defined so that for any iteration, every point in $\mathcal{X}$ is $\epsilon_{T}$-saturated. 
  Suppose $\epsilon \in (0,+\infty)$ is selected s.t.~$\epsilon_t$ is nondecreasing with $t$.
  Then for any $k \geq 1$, both Case (a) and Case (b) from Proposition~\ref{prop:eddp_next_sat} still hold w.r.t.~the accuracy tolerances $\epsilon_t$ defined in this proposition.
\end{proposition}

We also provide a numerical value to ensure every point is $\epsilon_{T}$-saturated. The proof is nearly identical to that of Lemma~\ref{lem:eddp_every_pt_is_sat}, so we omit it. 
\begin{lemma} 
  For any $k \geq 0$ and $x \in \mathcal{X}$, then $\oV^k(x) - \uV^k(x) \leq (1-\lambda)^{-1} [\overline{h}-\underline{h}]$.
\end{lemma}

Finally, observe that Lemmas~\ref{lem:eddp_find_or_converge} and~\ref{lem:eddp.1}, which are used for the convergence of the CE-Inf-EDDP (Algorithm~\ref{alg:alg5}) in Theorem~\ref{thm:thm2b} still hold for Gap-Inf-EDDP. 
Therefore, we have similar convergence guarantees.
\begin{theorem} \label{eq:tqgp_a1}
  Let $\overline{M} := \lambda \overline{M}_0 + M_h$.
  For any $T \geq 1$ and  $\epsilon \in \big(0, \frac{\lambda^T(\overline{h}-\underline{h})}{\overline M + \underline M} \big)$, 
  Gap-Inf-EDDP returns a feasible solution $x^k_0$ such that
  \begin{talign*} 
    F(x^k_0) - F^* 
    \leq
    \epsilon_0
    \leq
    \frac{(\overline{M}+\underline{M})\epsilon}{1-\lambda} + \frac{\lambda^{T}(\overline{h}-\underline{h})}{1-\lambda},
  \end{talign*}
  in at most $K$ iterations, where $K := 4T(D/\epsilon+1)^{n}$.
\end{theorem}
This result suggests that using the upper-bound model provides the same worst-case iteration complexity as with just the cutting-plane model from Theorem~\ref{thm:thm2b}.
Indeed, our later numerical experiments will show both methods have similar iteration complexity.
The next proposed search selection strategy offers a simpler strategy with probabilistic guarantees.

\subsection{Infinite-horizon stochastic dual dynamic programming} \label{sec:sec3}
We introduce a randomized CE-Inf-EDDP, which is similar to the celebrated stochastic dual dynamic programming~\cite{pereira1991multi,shapiro2011analysis,philpott2008convergence}, or SDDP.
CE-Inf-SDDP (Algorithm~\ref{alg:alg6}) selects trial points uniformly at random over scenarios $1,\ldots,N$ instead of deterministically like with CE-Inf-EDDP.
Because the trial points are now randomly selected, the algorithm can be significantly simplified because we do not need the accuracy level $\epsilon$ as an input, nor do we need to store/update the table $\Seps$.

\begin{algorithm}[h]
\begin{algorithmic}[1]
\Procedure{CE-Inf-SDDP}{$T$}
    \State{Run Inf-EDDP (Algorithm~\ref{alg:alg2}) but replace Line~\ref{line:line8} with
      \begin{talign*} 
        x^k &= \begin{cases} x^k_0 &: k~\mathrm{mod}~(2T) = 1 \\ 
          x^k_i & : \text{o/w, where}~ i \stackrel{\text{unif}}{\sim} \{\revone{0},\ldots,N\} \end{cases},
      \end{talign*}\label{line:sddp_selection}
      \revone{and also remove Lines~\ref{line:line1},~\ref{line:line6} to~\ref{line:line7},~\ref{line:z_def}, and~\ref{line:line5}}.
      }
\EndProcedure
\end{algorithmic}
\caption{}
\label{alg:alg6}
\end{algorithm}

Let us briefly discuss the convergence proof of CE-Inf-SDDP, since it can be adapted from~\cite[Theorem 3]{lan2022complexity} into CE-Inf-EDDP.
At a high level, we want CE-Inf-SDDP to mimic CE-Inf-EDDP.
Case 1 is easier to show this because its analysis decomposes the iterations $[1,k]$ into $2T$-length sequences $[\ell \cdot (2T)+1, (\ell+1) \cdot (2T)]$ (c.f.~Lemma~\ref{lem:eddp_find_or_converge_first}). 
Then we only need to show the CE-Inf-SDDP mimics CE-Inf-EDDP within sequences, and the probability of success within a single sequence is exponentially small in $T$.
In contrast, Case 2 views the first $k$ iterations as a single sequence $[1,k]$, and the probability of success can be exponentially small in $k$.
We now state the complexity of CE-Inf-SDDP.
\begin{theorem} \label{thm:tqgp_a1}
  For any $T \geq 1$ and $\epsilon \in \big(0, \frac{\lambda^T(\overline{h}-\underline{h})}{M + \underline{M}} \big]$, CE-Inf-SDDP generates a feasible solution $x_0^k$ s.t.
  \begin{talign*} 
    F(x^k_0) - F^* 
    \leq 
    \frac{(M + \underline{M})\epsilon}{1-\lambda}
    +
    \frac{\lambda^{T}(\overline{h} - \underline{h})}{1-\lambda},
  \end{talign*}
  in at most $K$ iterations, where $K$ is a random variable with an expected value of $\mathbb{E}[ K ] \leq \bar{K}_\text{SDDP} := 4T \bar K N^{2T}$ and $\bar K := (D/\epsilon+1)^n$.
  In addition, for any $\alpha \geq 1$, we have
  \begin{talign*}
      \mathrm{Pr}\{K  \geq \alpha \cdot K_{\text{SDDP}} + 1 \} 
      \leq 
      \mathrm{exp}\big\{ - \frac{4(\alpha-1)^2T \bar K}{\alpha \cdot N^{2T}} \big\}.
  \end{talign*}
\end{theorem}
Similar to~\cite[Theorem 3]{lan2022complexity}, the complexity of CE-Inf-SDDP can be exponentially worse in $T$ compared to CE-Inf-EDDP.
However, numerical evidence suggests this bound may be over-conservative.
Now, a standing assumption made so far in all our algorithms is that the subproblems~\eqref{eq:eq5} are simple enough (e.g., a deterministic LP) to be solved exactly.
Next, we consider the case where the subproblem can only be solved inexactly.

\subsection{Convergence rates for inexact subproblem solutions}
Our goal here is to understand how CE-Inf-EDDP performs when the subproblems~\eqref{eq:policy} and~\eqref{eq:eq5} from the forward step can only be solved inexactly.
To do so, we need a notion of inexactness, where our definition will use the feasible region $X(\cdot,\xi)$ from~\eqref{eq:qgp_a01} and the samples $\tilde{\xi}_i = \{\tilde{\B{A}}_i, \tilde{\B{B}}_i, \tilde{\B{Q}}_i, \tilde{\B{b}}_i, \tilde{\B{c}}_i, \tilde{\B{r}}_i\}$ from~\eqref{eq:saa_xi} for the SAA problem.
\begin{definition} \label{def:saderr.1}
  The pair $(\bar{x}^k_i, \bar{y}^k_i)$ is an $\epsilon$-approximate solution to subproblem~\eqref{eq:eq5} for scenario $i = 1,\ldots,N$ and iteration $k \geq 1$ when
  \begin{enumerate}
    \item $\uF^k_i(\bar{x}^k_i, \tilde{c}_i) - \uF^k_i(x^{k,*}_i, \tilde{c}_i) \leq \epsilon$, where $x^{k,^*}_i$ is the optimal primal solution to~\eqref{eq:eq5};
    \item $\tilde{\B B}_i^T \bar{y}^k_i$ is an $\epsilon$-subgradient satisfying
      \begin{talign*}
        \tv^k_i(x) \geq \tv^k_i(x^{k-1}) + \langle \tilde{\B B}_i^T \bar{y}^k_i, x-x^{k-1} \rangle - \epsilon,~\ 
        \forall x \in \mathcal{X}(\bar \epsilon);
      \end{talign*}
    \item The constraints are $\epsilon$-satisfied, or that $\tilde{\B A}_i \bar{x}^k_i - \tilde{\B B}_ix^{k-1} - \tilde{\B b}_i + \delta \in \mathcal{K}$ and $\phi(\bar{x}^k_i, \tilde{r}_i) + \delta_\phi \leq \tilde{Q}_i x^{k-1}$ for some $\delta \in \mathbb{R}^m$ and $\delta_\phi \in \mathbb{R}^{m_\phi}$ such that $\|\delta \|_\infty + \|\delta_\phi\|_\infty \leq \epsilon$.
  \end{enumerate}
  \sloppy Similarly, we say the here-and-now solution $\bar{x}_0^k$ is an $\epsilon$-approximate primal solution to subproblem~\eqref{eq:policy} when items 1 and 3 above hold w.r.t~the first-stage fixed data, $\xi^1 = \{A^1,B^1,Q^1,b^1,c^1,r^1\}$.
\end{definition}
Since inexact solutions may not be feasible, we need to extend the feasible set to include approximately feasible solutions, defined as
\begin{talign*}
    X(u,\tilde \xi_i,\delta,\delta_\phi) = \{ x \in \bar{X} : \tilde{\B A}_ix - \tilde{\B B}_iu - \tilde{\B b}_i  + \delta \in \mathcal{K},~\phi(x,\tilde{r}_i) +\delta_\phi \leq \tilde{Q}_i u\},
\end{talign*}
where $\delta \in \mathbb{R}^m$ and $\delta_\phi \in \mathbb{R}^{\phi}$ are some arbitrary vectors. 
We also define the sets of bounded perturbations for some $\Delta > 0$, 
\begin{talign*}
  \mathcal{D}_i(u) &= \{ (\delta,\delta_\phi) : X(u, \tilde \xi_i,\delta,\delta_\phi) \ne \emptyset \} \\
  \mathcal{D}_i(u, \Delta) &= \mathcal{D}_i(u) \cap \big \{(\delta,\delta_\phi) : \| \delta \|_\infty + \|\delta_\phi\|_\infty \leq \Delta \big\}.
\end{talign*}
The first set considers all perturbations $(\delta,\delta_\phi)$ such that the set $X(u, \tilde \xi_i, \delta, \delta_\phi)$ is non-empty. 
The second set restricts these perturbations to have bounded norms. 
We make the following assumption to help carry out a sensitivity analysis, where we recall $\mathcal X(\bar \epsilon)$ from~\eqref{eq:qgp_a9} and recall the scalar $\bar \epsilon \in (0,+\infty)$ from Assumption~\ref{asmp:eddp_asmp_relint}
\begin{assumption} \label{asmp:eps_0}
   There exists $\Delta \in (0,+\infty)$ s.t.~for every $x \in \mathcal{X}(\bar \epsilon)$, $\emptyset \ne \mathcal{D}_i(x, \Delta) \subseteq \mathrm{rint} \mathcal{D}_i(x)$.
\end{assumption}
This is a regularity condition that is similar to Assumption~\ref{asmp:eddp_asmp_relint}.
This assumption can be satisfied in certain situations by enlarging the feasible set via perturbations. 
For example, if the feasible region contains only affine constraints, $X(u,\tilde \xi_i) = \{x \in \bar X~ \vert~ \tilde{\B A}_ix - \tilde{\B B}_iu \leq \tilde{\B b}_i \}$, which is non-empty for all $u \in \mathcal X(\bar \epsilon)$, then the assumption can be satisfied when replacing $\tilde{\B b}_i$ with $\tilde{\B b}_i + \Delta \cdot \mathbf{1}$.

We can now establish an important Lipschitz property.
The proof can be shown similarly to Lemma~\ref{lem:eddp_asmp_smooth1}, so we skip it.
\begin{lemma} \label{lem:pert_feas}
  For any $k \geq 1$ and $i=0,\ldots,N$, define the perturbed problem,
  \begin{talign} \label{eq:qgp_a10}
    g_i^k(\delta,\delta_\phi) := \min_{x \in X(x^{k-1}, \tilde \xi_i, \delta, \delta_\phi)} F_i^k(x), \quad  \forall~(\delta,\delta_\phi) \in \mathcal{D}_i(x^{k-1}, \Delta).
  \end{talign}
  There exists a universal constant $M_D < +\infty$
  s.t.
  \begin{talign*}
    \vert g_i^k(\delta,\delta_\phi) - g_i^k(\delta',\delta_\phi') \vert 
    &\leq 
    M_D [\|\delta - \delta'\|_\infty + \|\delta_\phi - \delta'_\phi\|_\infty],
    \\
    &\forall (\delta,\delta_\phi), (\delta',\delta_\phi') \in \mathcal{D}_i(x^{k-1}, \Delta).
  \end{talign*}
\end{lemma}
Note that $g_i^k(\cdot,\cdot)$ will only be used for the analysis and need not be computed.
A finite bound for $M_D$ can be established, similarly to the bound for $M_V$ from Lemma~\ref{lem:eddp_asmp_smooth1}. 
We omit this derivation for simplicity, since it requires updating the feasible region $\mathcal{X}(\bar \epsilon)$ in~\eqref{eq:qgp_a9} to hold w.r.t.~the enlarged feasible region $X(\cdot,\tilde{\xi}_i,\delta,\delta_\phi)$.

We begin our convergence analysis by relating the value function $V$ with its cutting-plane model, $\uV^k$. It is important to note that, unlike the setting when we solve the subproblems exactly (see Lemma~\ref{lem:eddp_relation}), the cutting-plane model $\uV^k$ with inexact supporting hyperplanes is no longer an underestimation of $V$. 
We skip the proof since it is a straightforward adaptation of Lemma~\ref{lem:eddp_relation} with the addition of the subgradient error from Definition~\ref{def:saderr.1}.
\begin{lemma} \label{lem:sto_eddp_relation}
  Suppose that at every iteration, an $\epsilon$-approximate solution is found for subproblem~\eqref{eq:eq5}, and an $\epsilon$-approximate primal solution is found for~\eqref{eq:policy}.
  Then for any iteration $k \geq 1$ and $i=0,\ldots,N$, 
 \begin{talign*}
   \uV^{k-1}(x) \leq \uV^k(x) 
   \leq 
   \frac{1}{N} \sum\limits_{i=1}^N \tv_i^k(x) + \epsilon 
   \leq 
   V(x) + \frac{1-\lambda^k}{1-\lambda} \epsilon,
   &\ \forall x \in \mathcal{X}(\bar{\epsilon}) \\ 
   \uF^{k-1}_i(x) \leq \uF^k_i(x) 
   \leq F_i(x) + \frac{\lambda(1-\lambda^k)}{1-\lambda} \epsilon, 
   &\ \forall x \in \mathcal{X}(\bar{\epsilon}), 
 \end{talign*}
 where $F_0 := F$ is from~\eqref{eq:qgp_a3} and $\uF_0^k := \uF^k$ is from~\eqref{eq:policy}.
\end{lemma}

The next result is nearly identical to Proposition~\ref{prop:eddp_next_sat}, where we provide sufficient conditions for a search point $x^k$ to become saturated. 
The only difference between the result below and Proposition~\ref{prop:eddp_next_sat} is that the saturation levels $\epsilon_t$ are modified to account for primal and dual errors.
\begin{proposition} \label{prop:sto_eddp_next_sat}
  Let $\epsilon_t$ follow the recursive definition $\epsilon_{t-1} = (M + \underline M + M_D + 1)\epsilon  + \lambda \epsilon_t$, and let $\epsilon_{T}$ be defined so that for any iteration, every point in $\mathcal{X}$ is $\epsilon_{T}$-saturated. 
  Suppose $\epsilon \in (0,+\infty)$ is selected s.t.~$\epsilon_t$ is nondecreasing with $t$, and that at every iteration, an $\epsilon$-approximate solution is found for subproblem~\eqref{eq:eq5}.
  Then for any $k \geq 1$, the following two conditions hold:
  \begin{enumerate}[(a)]
      \item \textit{Case (a)} For any $x \in \mX$, if $\Seps^k (\revone{\Geps}(\revtwo{x})) = t$, then
           \begin{talign*} 
             \lambda [ V(\revtwo{x}) - \uV^{k-1} (\revtwo{x}) ] 
             =
             F_i(\revtwo{x}) - \uF_i^k (\revtwo{x}) 
            \leq
            (M + \underline M) \epsilon + \lambda\epsilon_t,
           \end{talign*}
           where $F_0 := F$ is from~\eqref{eq:qgp_a3} and $\uF_0^k := \uF^k$ is from~\eqref{eq:policy}.
           \sloppy \revtwo{Furthermore, when $\max_{i =1,\ldots,N} \Seps^k (\revone{\Geps}(x^k_i)) = t$, where the trial point $x^k_i$ is from Line~\ref{line:line2}, then $x^{k-1}$ is $\epsilon_{t-1}$-saturated at iteration $k$, or equivalently,
           \begin{talign*} 
              V(x^{k-1}) - \uV^k(x^{k-1}) \leq \epsilon_{t-1}.
           \end{talign*}}
        \item \textit{Case (b)} The same Case (b) as Proposition~\ref{prop:eddp_next_sat}, which we repeat for completeness. 
        For any $x \in \mathcal X$, if $\Seps^{k+1}(\revone{\Geps}(x)) = t < T$, then there exists a previous trial point $x^j$ (where $j \leq k$) s.t.~$\|x-x^j\|_\infty \leq \epsilon$ and $x^j$ is $\epsilon_{t'}$-saturated for any $t' \in [t,T]$ at iteration $k+1$.
  \end{enumerate}
\end{proposition}
\begin{proof}
  The proof is nearly identical to Proposition~\ref{prop:eddp_next_sat}. 
  We only need to modify~\eqref{eq:qgp_a43} and~\eqref{eq:qgp_a42}, which are bounds on $V(x^{k-1})$ and $\uV^k(x^{k-1})$, respectively.  

  To that end, let $(x^k_i,y^k_i)$ be the computed $\epsilon$-approximate solution to~\eqref{eq:eq5}, and let $\delta$ be the vector where $x^k_i \in X(x^{k-1}, \tilde \xi_i, \delta, \delta_\phi)$.
  We introduce $\bar{x}^k_i$ as the optimal solution to $g_i^k(\mathbf 0, \mathbf{0})$ and $\tilde{x}^k_i$ as the optimal solution to $g_i^k(\delta, \delta_\phi)$ from~\eqref{eq:qgp_a10}. 
  Then 
  \begin{talign*}
    V(x^{k-1})
    &\leq
    \frac{1}{N} \sum_{i=1}^N [F_i(x^k_i) + ( F_i(\xeps) - F_i(x^k_i)) ] \\
    &\leq
    \frac{1}{N} \sum_{i=1}^N [F_i(x^k_i) + ( F_i(\xeps) - F_i(\tilde{x}^k_i)) ] \\
    &=
    \frac{1}{N} \sum_{i=1}^N [F_i(x^k_i) + ( g_i^k(\mathbf{0}, \mathbf{0}) - g_i^k(\delta, \delta_\phi) )]
    \\ &\leq
    \frac{1}{N} \sum_{i=1}^N  F_i(x^k_i) + M_D \epsilon,
  \end{talign*}
  where the first inequality is by~\eqref{eq:qgp_a43}, 
  the second inequality by observing $x^k_i$ is a feasible solution to the optimization problem given by $g_i^k(\delta,\delta_\phi)$ because $x_i^k \in X(x^{k-1}, \tilde \xi_i, \delta, \delta_\phi)$,
  the equality by recalling $g_i^k(\cdot, \cdot)$ from~\eqref{eq:qgp_a10} and our construction of $\bar{x}^k_i$ and $\tilde x^k_i$, 
  and the last inequality is by Lemma~\ref{lem:pert_feas} and $\|\delta\|_\infty + \|\delta_\phi \|_\infty \leq \epsilon$ by the fact that $(x^k_i,y^k_i)$ is an $\epsilon$-approximate solution.

  Now, define $\hat x^k_i$ as the optimal solution to $\tv_i^k(x^{k-1})$, where the latter is from~\eqref{eq:eq5}.
  Then
  \begin{talign*}
    \uV^k(x^{k-1})
    &\geq 
    \frac{1}{N} \sum_{i=1}^N \tv_i^k(x^{k-1}) \\
    &= 
    \frac{1}{N} \sum_{i=1}^N [\uF_i^k(x^k_i) + ( \uF_i^k(\hat x^k_i) - \uF_i^k(x^k_i) )]
    \geq 
    \frac{1}{N} \sum_{i=1}^N \uF_i^k(x^k_i) - \epsilon,
  \end{talign*}
  where the first inequality is by~\eqref{eq:qgp_a42}, the equality is by definition of $\hat x^k_i$ and $\tv_i^k$, and the last inequality is by $x^k_i$ being an $\epsilon$-approximate solution.
\end{proof}

We are now ready to establish the convergence of CE-Inf-EDDP with inexact subproblems.
\begin{theorem} \label{thm:thm3}
  For any $T \geq 1$ and $\epsilon \in \big(0, \frac{\lambda^T(\overline{h}-\underline{h})}{M + \underline M + M_D + 1} \big]$, 
  let $\epsilon_t$ be \revone{recursively defined as in Proposition~\ref{prop:sto_eddp_next_sat}}.
    Suppose that at every iteration, an $\epsilon$-approximate solution is found for subproblem~\eqref{eq:eq5} and an $\epsilon$-approximate primal solution is found for~\eqref{eq:policy}.
    Then the CE-Inf-EDDP returns a solution ${x}^k_0$ s.t.
    \begin{talign*} 
        F(x^k_0) - F^* \nonumber \leq \epsilon_0 + \frac{\epsilon}{1-\lambda}
        & \leq 
        \frac{(M + \underline{M} + M_D + 2)\epsilon }{1-\lambda} 
        + \frac{\lambda^{T}[\overline{h}-\underline{h}]}{1-\lambda}
    \end{talign*}
    in at most $K$ iterations, where $K:= 4T(D/\epsilon+1)^n$.
\end{theorem}
\begin{proof}
  The proof is mostly the same as Theorem~\ref{thm:thm2b}, except with two main differences.
  First, $\epsilon_0 \leq \frac{(M + \underline{M} + M_D + 1) \epsilon }{1-\lambda} + \frac{\lambda^T[\overline{h}-\underline{h}]}{1-\lambda}$,
  where we use the new recursive definition of $\epsilon_t$. Second, we modify~\eqref{eq:F_optimality_relationship} with
  \begin{talign*}
      F(x^k_0) - F(x^*)
      &\leq
      F(x^k_0) - \uF^k_0(x^*) + (1-\lambda)^{-1}\epsilon \\
      &\leq
      F(x^k_0) - \uF^k_0(x^k_0) + (1-\lambda)^{-1}\epsilon
      \leq
      \epsilon_0 + (1-\lambda)^{-1}\epsilon,
  \end{talign*}
  where we used Lemma~\ref{lem:sto_eddp_relation}.
\end{proof}
Compared to CE-Inf-EDDP with exact sub-problem solutions from Theorem~\ref{thm:thm2b}, Theorem~\ref{thm:thm3} above says inexact solutions only increase the final sub-optimality gap additively by $\frac{(1+M_D)\epsilon}{1-\lambda}$, which scales linearly with the sub-problem error $\epsilon$.
Furthermore, convergence guarantees with inexact solutions can be similarly extended to Gap-Inf-EDDP and CE-Inf-SDDP as well.
Next, we examine the numerical performance of the methods studied so far.

\section{Numerical Experiments} \label{sec:experiments}
\sloppy We now apply our various infinite-horizon dual dynamic programming (Inf-DDP) algorithms -- Inf-EDDP, CE-Inf-EDDP, Gap-Inf-EDDP, and CE-Inf-SDDP -- to solve the inventory (i.e.,~newsvendor) and hydrothermal planning problem.
We also compare them to finite-horizon methods (with horizon length $T$ specified later) of EDDP and cyclic-SDDP (Cyc-SDDP), the latter of which runs SDDP with a single cutting-plane model.
Cyc-SDDP is similar to running SDDP on the so-called cyclic policy graphs~\cite{dowson2020policy} for modeling stationary infinite-horizon problems.
The problems are set up similarly to~\cite{shapiro2020periodical}, where the number of scenarios is $N=50$, and we use a small $\lambda=0.8$ and a large discount factor $\lambda=0.9906$.
\revone{See~\cite{shapiro2023dual} for further explanations behind the discount factor choices.}
All methods compute a lower-bound estimate directly from the optimal value to the here-and-now problem~\eqref{eq:policy} as well as an upper-bound estimate. 
\revone{To estimate an upper bound of the optimal value, we take the aforementioned lower bound and add to it $\oV^k(x^k_0) - \uV^k(x^k_0)$, where we recall that $x^k_0$ is the optimal solution to~\eqref{eq:policy}.
In other words, we simply replace the lower bound estimate $\uV^k(x^k_0)$ with the upper bound estimate $\oV^k(x^k_0)$ defined in~\eqref{eq:qgp_a15}}.

All subproblems are solved using Gurobi v12.0.3 on one core of an Intel Xeon 6972P processor with 8GB of memory.
Code and more details can be found at \url{https://github.com/jucaleb4/hddp}.

\subsection{An infinite-horizon inventory problem}
In the newsvendor or inventory planning problem~\cite{shapiro2021lectures}, the goal is to control the inventory level over a multi-period horizon,
where the cost function is
\begin{talign} \label{eq:inv_cost}
  h(x^t,\tilde{\xi}_i) := \langle c, x^t-y^t \rangle + \langle b, [-y^t]_+ \rangle + \langle h, [y^t]_+ \rangle,
\end{talign}
where $[\cdot]_+ := \max\{0,\cdot\}$, $\tilde{D}^t := D(\tilde \xi_i) \in \mathbb{R}^n$ is the random demand at time $t$. 
There are two variables: $y^t \in \mathbb{R}^n$ is the inventory at the start of time $t$, which is subject to (s.t.) $y^t = x^{t-1} - \tilde D^t$, and $x^t \in \mathbb{R}^n$ is the inventory after replenishment at time $t$ s.t.~$x^t \geq y^t$.
We also include box constraints for these variables.
Vectors $c,b,h$ are the ordering, backlog, and holding costs, respectively, of appropriate size.
We first consider the simple $n=1$ case, and we ran with both a small discount factor $\lambda = 0.8$ (and $T=24$) and a large factor $\lambda = 0.9906$ (and $T=120$).

\begin{figure}[H]
  \begin{subfigure}{.49\textwidth}
      \includegraphics[width=0.8\linewidth]{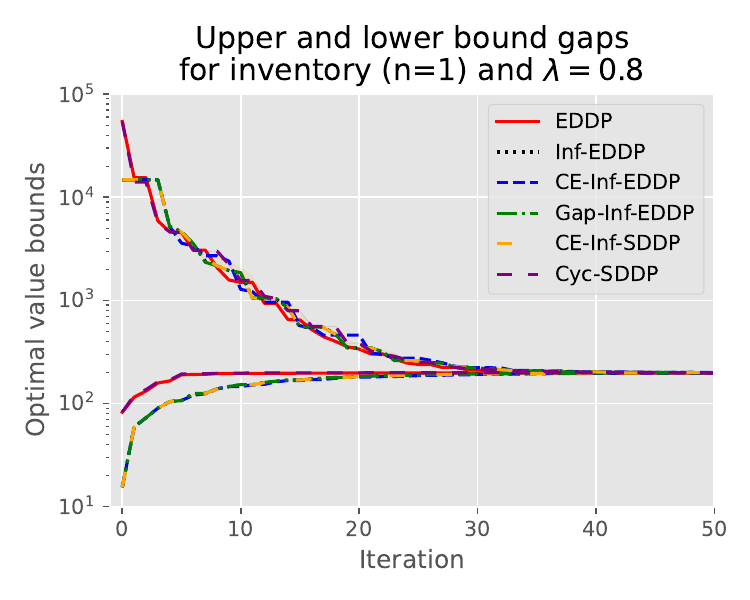}
      \centering
   \end{subfigure}%
  \begin{subfigure}{.49\textwidth}
    \includegraphics[width=0.8\linewidth]{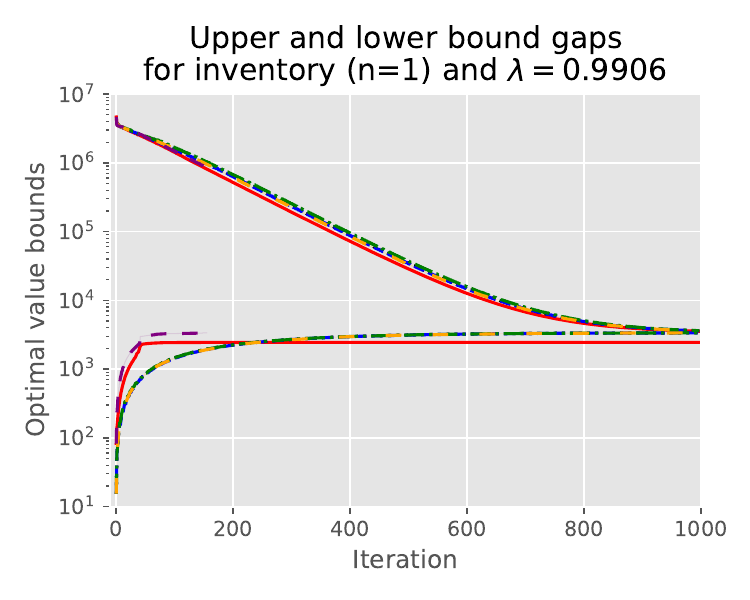}
    \centering
  \end{subfigure}
  \caption{Upper and lower bound gap estimates of the optimal value for a 1-D inventory problem for different DDP algorithms.
  CE-Inf-SDDP and Cyc-SDDP are repeated over 10 trials, where the line represents the mean and the shaded region represents the 95\% confidence interval.
  \revone{Some methods have near-identical performances, and thus their lines may overlap. 
  See the text for discussions on similar performing methods.}
  }
  \label{fig:newsvendor_gaps}
\end{figure}

As seen in Figure~\ref{fig:newsvendor_gaps}, all methods can be solved to near optimality.
For $\lambda=0.8$, all methods obtain a final relative optimality gap, defined as $(\text{upper bound} - \text{lower bound})/\text{lower bound}$, of at most 5.2e-03.
In contrast, for $\lambda=0.9906$, \revone{our family of Inf-EDDP has the tightest final relative optimality gap of 7.0e-2.
While EDDP and Cyc-SDDP show promising early performance, their final relative optimality gaps are 3.8e-1 and 2.4e+2, respectively, which cannot match our methods in the long term.}
In particular, Cyc-SDDP has the largest gap because it terminates before iteration 200 due to exceeding the two-hour time limit.
It seems our methods explore the feasible region more adaptively, yielding smaller optimality gaps in the long term.

\begin{figure}[h!]
  \begin{subfigure}{.32\textwidth}
      \includegraphics[width=\linewidth]{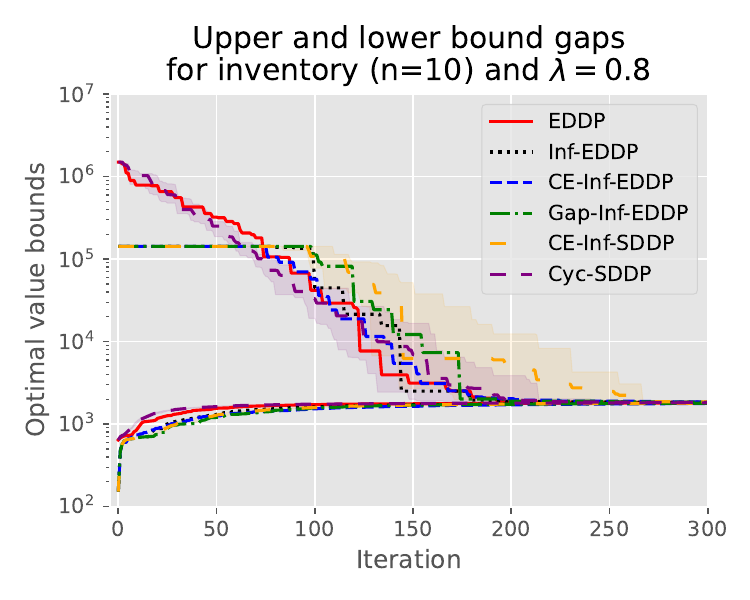}
      \centering
   \end{subfigure}
  \begin{subfigure}{.32\textwidth}
    \includegraphics[width=\linewidth]{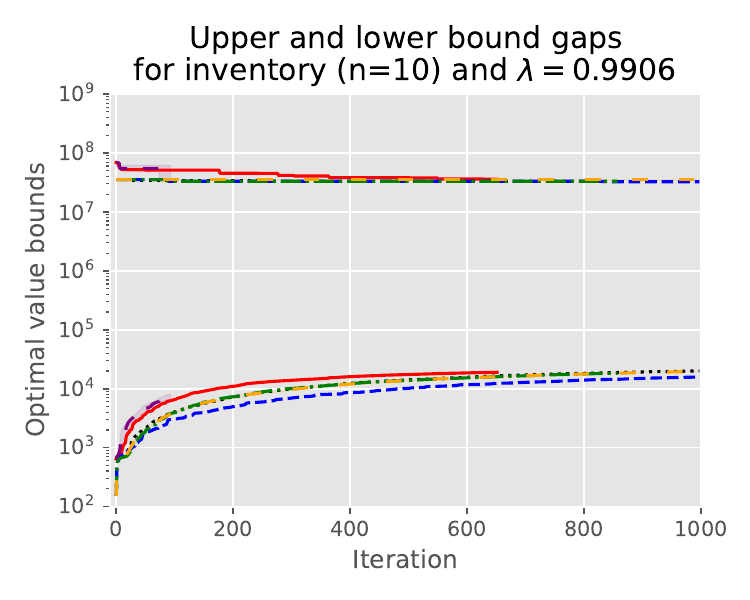}
    \centering
  \end{subfigure}
  \begin{subfigure}{.32\textwidth}
    \includegraphics[width=\linewidth]{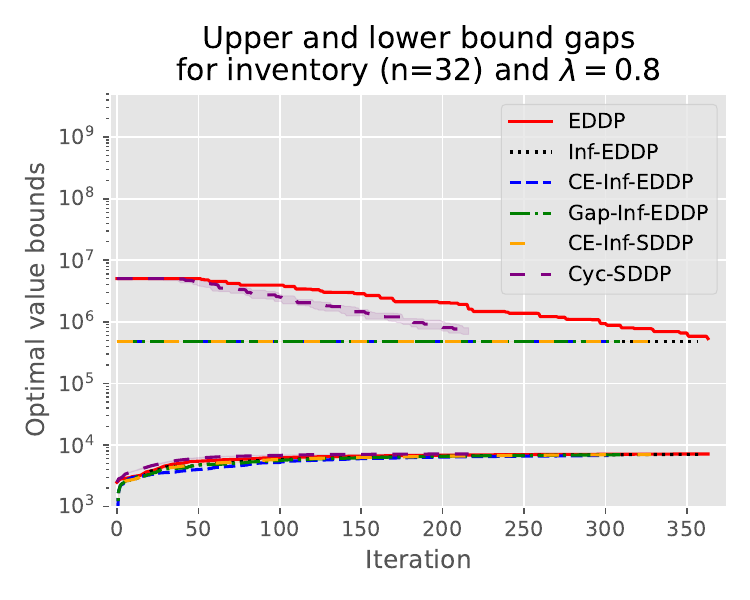}
    \centering
  \end{subfigure}
  \caption{Upper and lower bound gap estimates of the optimal value for an n-D inventory problem (with $n =10,32$) for different DDP algorithms. Same setup as Figure~\ref{fig:newsvendor_gaps}.}
  \label{fig:bigger_newsvendor_gaps}
\end{figure}

\revone{For larger instances with more inventory $n$, the time to get to near certifiable optimality or the final relative optimality gap gets larger, as displayed in Figure~\ref{fig:bigger_newsvendor_gaps}.
In the large instance of $n=10$ and $\lambda = 0.8$, all methods can get a final relative optimality gap of at most 4.32e-02, but the number of iterations to do is about 6x more than in $n=1$. 
For $n=10$ and $\lambda = 0.9906$ or $n=32$ and $\lambda = 0.8$, no methods can solve close to optimality, with the smallest final optimality gap being 1.84e+03 and 9.11e+01, respectively.}

\begin{table}
   \centering
   \begin{tabular}{@{}lllrrrr@{}} \toprule
    Gap & $n$ & $\lambda$ & CE-Inf-EDDP & EDDP & Cyc-SDDP & Gap-Inf-EDDP \\ \midrule
    \tikzxmark & 1 & 0.8 & \textbf{8.19e-3} & 4.30e-1 & 2.92e+0 & 2.11e+0 \\
    \tikzxmark & 1 & 0.9906 & \textbf{1.52e-2} & 4.53e+0 & 1.73e+1 & 2.48e+0  \\
    \tikzxmark & 10 & 0.8 & \textbf{3.06e-2} & 2.00e+0 & 9.72e+0 & 3.37e+0  \\
    \tikzxmark & 10 & 0.9906 & \textbf{3.96e-2} & 8.67e+0 & 7.28e+1 & 4.99e+0  \\ \midrule
    \tikzcmark & 1 & 0.8 & \textbf{1.37e-1} & 6.31e-1 & 3.13e+0 & 2.24e+0  \\
    \tikzcmark & 1 & 0.9906 & \textbf{4.99e-1} & 5.16e+0 & 1.74e+1 & 2.92e+0  \\
    \tikzcmark & 10 & 0.8 & \textbf{1.48e+0} & 7.48e+0 & 1.05e+1 & 4.81e+0  \\
    \tikzcmark & 10 & 0.9906 & \textbf{3.59e+0} & 1.10e+1 & 7.32e+1 & 8.39e+0 \\
    \bottomrule
   \end{tabular}
   \caption{Average runtime in seconds per iteration for solving the inventory problem with different sizes of $n$.
   The least time is bolded.
   Column ``Gap'' indicates whether the runtime for computing the optimality gap is included.
   We omit Inf-EDDP and CE-Inf-SDDP for simplicity, since their recorded runtimes are similar to CE-Inf-EDDP.}
   \label{tab:newsvendor_average_runtime}
\end{table}

Our family of Inf-EDDP methods so far has achieved competitive performance with existing DDP methods when comparing the optimality gap over the iterations.
It turns out the computational cost of our methods is much lower, as seen in the average runtime per iteration in Table~\ref{tab:newsvendor_average_runtime}.
Our CE-Inf-EDDP requires about 50x and 250x less runtime per iteration compared to EDDP for $\lambda = 0.8$ and $0.9906$, respectively, when excluding the optimality gap computation.
This matches our expected speedup on the order of the effective horizon length $T$ (recall $T=24$ and $T=120$ for $\lambda=0.8$ and $\lambda=0.9906$, respectively), because EDDP requires $T$ forward steps and $T-1$ backwards steps per iteration, whereas our methods only have a single forward and backward step.
We omitted CE-Inf-SDDP since it has nearly identical runtimes as CE-Inf-EDDP.
In particular, CE-Inf-SDDP achieved the largest relative speedup of 4.9\% on the instance $n=1$ and $\lambda=0.8$, with a runtime of 7.81e-3s without gap computations.
This confirms that the saturation table $\Seps$ to select search points in CE-Inf-EDDP has minimal computational overhead compared to the random sampling used in CE-Inf-SDDP.
Meanwhile, Cyc-SDDP has the largest runtime when gap computations are excluded. 
This is because it adds multiple cuts (or shares cuts) to the same cutting-plane model each iteration, resulting in larger sub-problems that require more solve time. 
In contrast, both EDDP and our methods only add one cut to any cutting-plane model each iteration.

Table~\ref{tab:newsvendor_average_runtime} also includes the average runtime when optimality gap computations are included. 
Here, the improvement is more modest, with about a 3x to 5x reduction in runtime for our methods. 
\revone{This large overhead comes from updating and evaluating the upper-bound model $\uV^k$, which helps compute the optimality gap.
At iteration $k$, updating $\oV^k(\cdot)$ by adding a constraint involves solving $N$ linear programs (LPs)~\eqref{eq:qgp_a14}, each with $O(N \cdot (k+n))$ variables and $O(N \cdot n)$ constraints, where recall $N=50$ is the number of scenarios, and $n$ is the dimension of the decision variable.
Evaluating $\oV^k(x)$ for a single $x$ involves solving $N$ LPs~\eqref{eq:upper_V}, each with $O(n)$ variables and $O(n+k)$ constraints.
In contrast, solving the sub-problems in Line~\ref{line:line2} in the forward step of a DDP method requires solving $N$ LPs (one for each scenario $i=1,\ldots,N$), each with $O(n)$ variables and $O(k)$ constraints, which is easier to solve than updating and computing $\uV^k(x^k_0)$ for the here-and-now solution $x^k_0$.
This also explains the longer runtime for Gap-Inf-EDDP.
Since it updates the saturation table $\Seps$ by evaluating $\oV^k(x^k_i)$ for each of the trial points $\{x^k_i\}_{i=0}^N$ in the added~\eqref{eq:qgp_a46}, the total number of LPs solved is $O(N^2)$, each with $O(n)$ variables and $O(n+k)$ constraints.}

\begin{figure}
  \begin{subfigure}{.49\textwidth}
    \includegraphics[width=\linewidth]{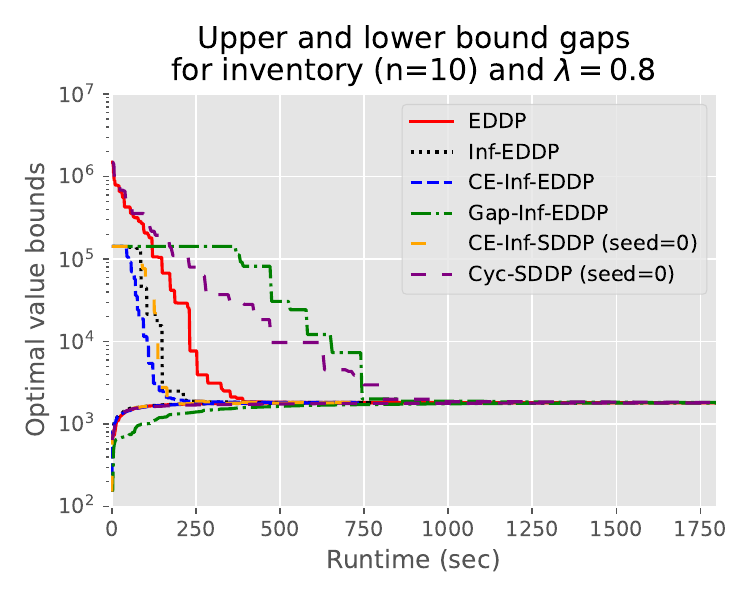}
    \centering
    \end{subfigure}
    \begin{subfigure}{.49\textwidth}
      \includegraphics[width=\linewidth]{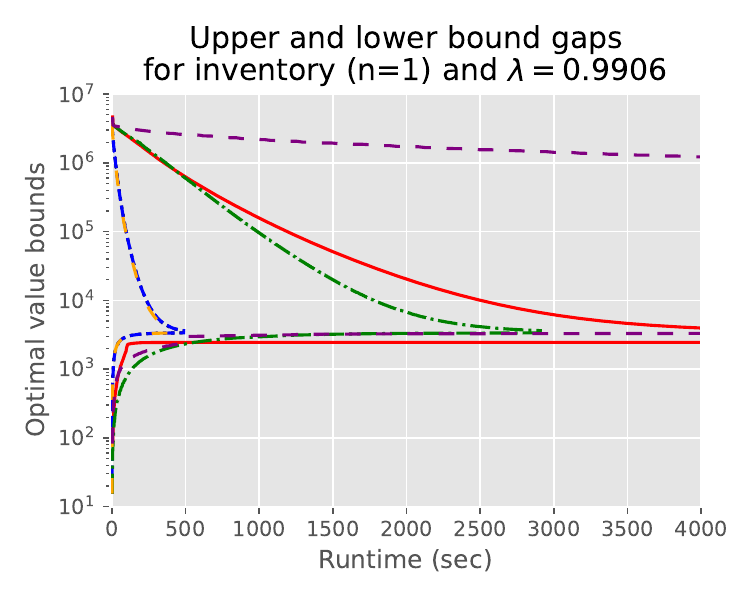}
      \centering
   \end{subfigure}%
  \caption{Upper and lower bound gap estimates of the optimal value w.r.t~runtime (including gap computations) for an inventory problem. Same setup as~\ref{fig:newsvendor_gaps}.}
  \label{fig:newsvendor_runtimegaps}
\end{figure}

To better understand the role of runtime, we added Figure~\ref{fig:newsvendor_runtimegaps}, which shows the optimality gap compared to runtime for two inventory instances.
In the instance of size $n=1$ and $\lambda = 0.9906$, our methods get a nearly tight optimality gap after about 500s of runtime, while it takes EDDP about 4000s.
In the instance $n=10$ and $\lambda=0.8$, CE-Inf-EDDP has tighter upper bounds than Inf-EDDP in the first 250 iterations, which showcases the effectiveness of the longer and adaptive exploration developed in subsection~\ref{sec:better_T}.
Note that these runtimes include gap computations computed every iteration.
Faster runtimes can be attained by periodically computing the gap or through parallel processing, since the LPs can be solved independently.

\subsection{An infinite-horizon risk-averse inventory problem}
We consider the risk-averse inventory problem, where the decision-maker prefers smaller inventory mismatch, \revone{which we model by a quadratic functional constraint.
This environment will test the ability of our methods to generalize towards nonlinear functional constraints.}
Our formulation uses the same cost and constraints from~\eqref{eq:inv_cost} with an additional quadratic constraint $(y^t)^2 - p^t \leq \tau$ for some tolerance $\tau > 0$, where $p^t \geq 0$ is a penalty variable with cost $C \cdot p^t$ to ensure feasibility (we set $C=100$). 
In our experiments, we set $\tau = 5$ since the policies from the risk-neutral setting prefer inventory levels of $y^t \approx 20$. 

\begin{figure}
  \begin{subfigure}{.49\textwidth}
      \includegraphics[width=\linewidth]{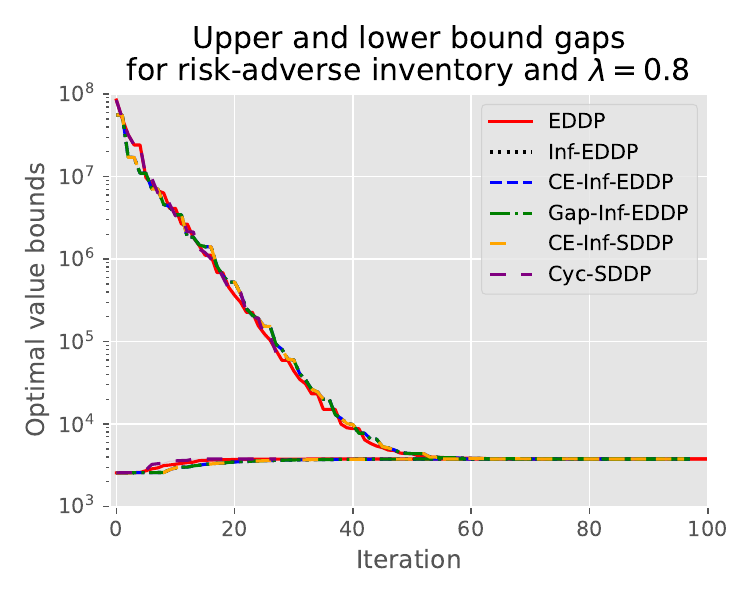}
      \end{subfigure}
    \begin{subfigure}{.49\textwidth}
      \includegraphics[width=\linewidth]{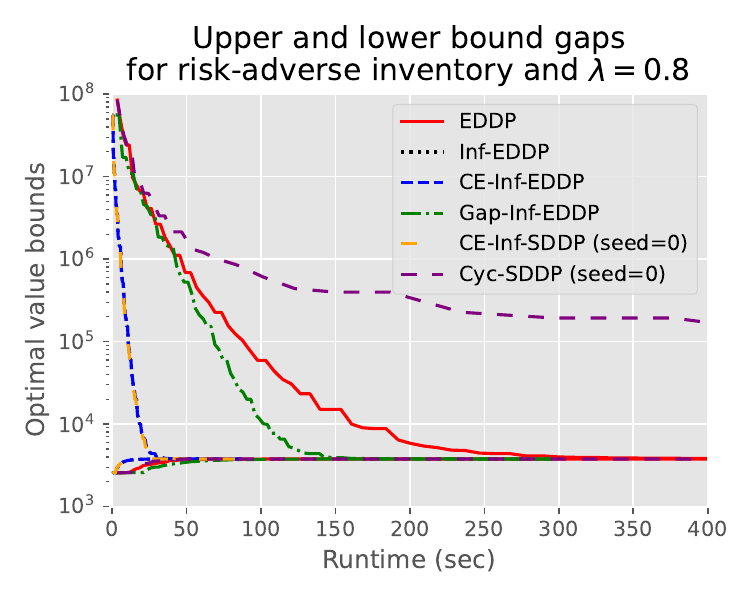}
      \centering
    \end{subfigure}%
    \caption{Upper and lower bound gap estimates of the optimal value w.r.t~both iterations and runtime (including gap computations) for a risk-averse inventory problem. Same setup as Figure~\ref{fig:newsvendor_gaps}.}
    \label{fig:ra_newsvendor}
\end{figure}
We plot similar lower- and upper-bound gap progressions in Figure~\ref{fig:ra_newsvendor}. 
Only $\lambda=0.8$ is shown since the runtime for $\lambda=0.9906$ with a longer time horizon is much longer to solve.
All methods except for Cyc-SDDP (since it can only achieve a relative optimality gap of 12.0 after the runtime limit of 2hrs) can achieve a tight optimality gap.
In particular, our family of Inf-EDDP methods has a relative optimality gap of at most 6.95e-5, while Inf-EDDP has a marginally larger relative gap of 1.92e-3.

\subsection{An infinite-horizon hydrothermal planning problem}
We now consider the infinite-horizon hydrothermal planning problem. 
\revone{This decision-making problem involves $n=4$ state variables for energy stored in hydro-reservoirs and $57$ control variables (i.e.,~decision variables that are not passed between stages) that correspond to previous energy stored, spillage, hydro-energy, energy deficit, thermal, and exchange variables.
Besides simple variable bound constraints, this problem has 13 equality constraints for flow and balance equations}.
See~\cite{shapiro2020periodical} for more details.

\begin{figure}[h!]
  \begin{subfigure}{.45\textwidth}
      \includegraphics[width=\linewidth]{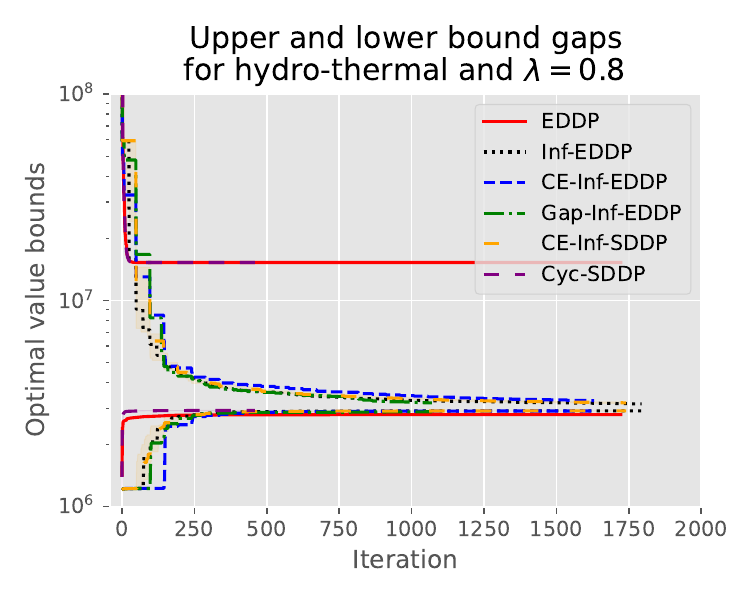}
      \centering
  \end{subfigure}
  \begin{subfigure}{.45\textwidth}
      \includegraphics[width=\linewidth]{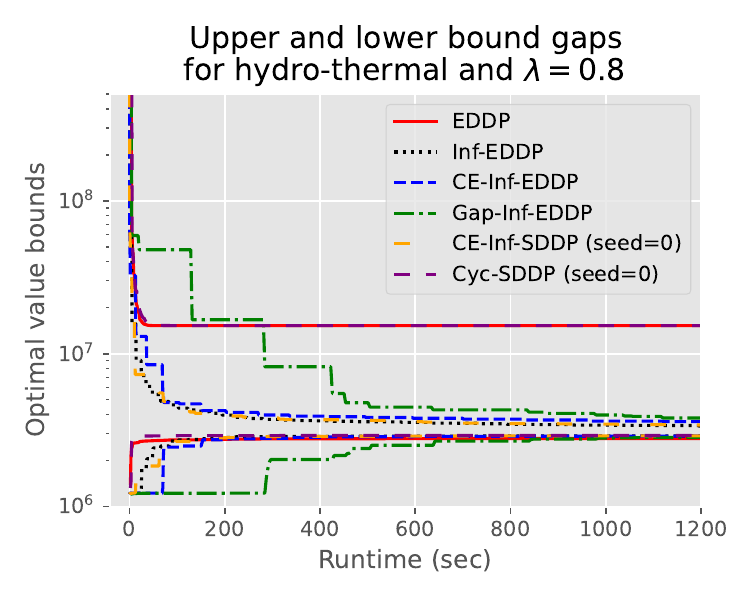}
      \centering
   \end{subfigure}%
   \newline
  \begin{subfigure}{.45\textwidth}
      \includegraphics[width=\linewidth]{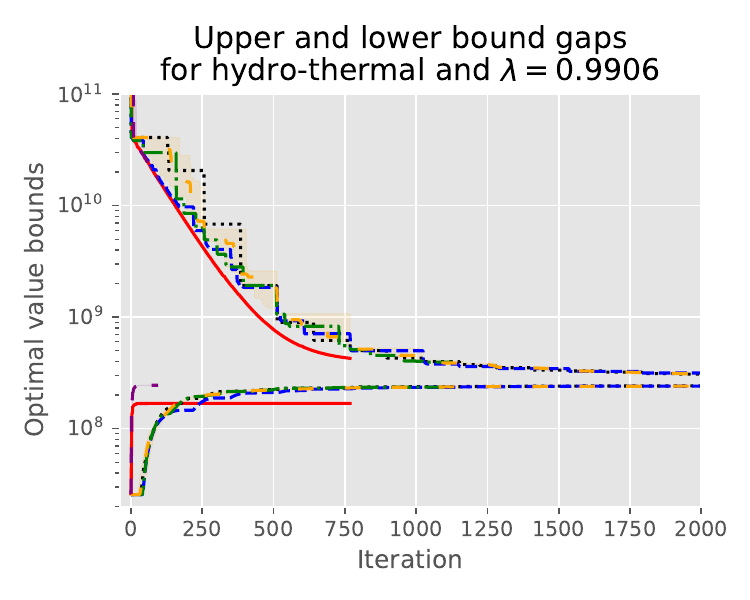}
      \centering
  \end{subfigure}
  \begin{subfigure}{.45\textwidth}
      \includegraphics[width=\linewidth]{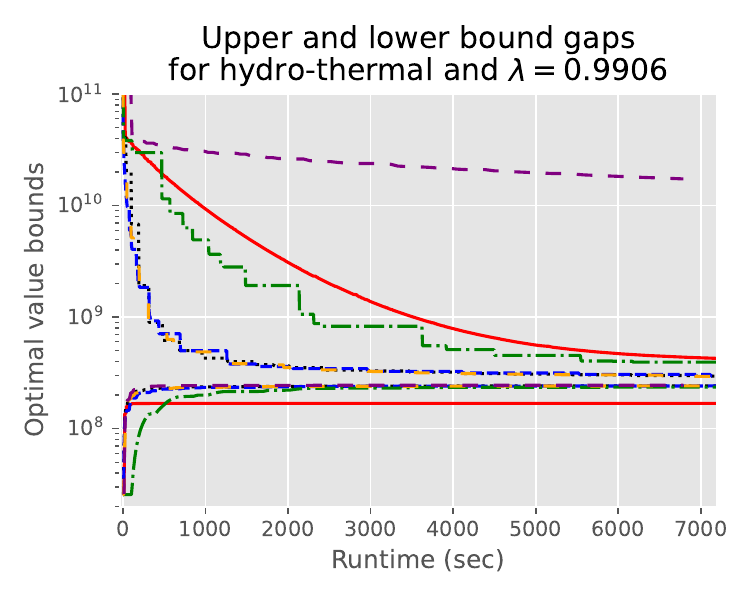}
      \centering
   \end{subfigure}%
   \caption{Optimality gaps for solving an infinite-horizon hydrothermal problem.
   Same setup as Figure~\ref{fig:newsvendor_gaps}.
   The left pair of plots shows the gap vs.~iteration count, while the right pair shows the gap vs.~runtime (including gap computations).}
\label{fig:hydro_gaps}
\end{figure}

The progress of the optimality gap over the different iterations is shown in Figure~\ref{fig:hydro_gaps}.
Our family of Inf-EDDP-type methods has a largest final relative optimality gaps of 5.18e-1 and 3.89e-1 for $\lambda=0.8$ and $\lambda=0.9906$, respectively, whereas Cyc-SDDP has about an order of magnitude larger final relative optimality gap of 4.22e+0 and 6.70e+0, respectively (EDDP yields similar gaps).
Similar to the inventory problem, the optimality gap for our family of Inf-EDDP methods converges more efficiently when measuring runtime instead, since each iteration of our Inf-EDDP-type methods has lower computational costs.
Likewise, Gap-Inf-EDDP can also attain similarly tight bounds, but its runtime is longer than CE-Inf-EDDP due to evaluating the upper bound more often.
In the instance with $\lambda=0.9906$, the longer per-iteration runtime of EDDP and Cyc-SDDP also explains why they reach the time limit and terminate early before iteration 800.
\revone{Note that in the instance with $\lambda=0.8$, EDDP and Cyc-SDDP's optimality gap stagnates early on.
This worse performance is because EDDP and Cyc-SDDP's search points $x^k$ may not change for larger iterations $k$.
More precisely, 131 out of the first 200 pairs of consecutive search points, $x^k$ and $x^{k+1}$, are identical for EDDP, while only 5 out of the first 200 pairs are identical for CE-Inf-EDDP.}
By exploring the feasible region adaptively and for longer, CE-Inf-EDDP can explore more of the feasible region to obtain smaller optimality gaps.

To understand the impact of a tight lower-bound model and/or optimality gap, we evaluate the in- and out-of-sample performance for a select number of methods. 
Both analyses are performed using a statistical upper-bound estimation over a finite horizon $T' = 1280$, where we perform $T'$ forward steps (i.e.,~solve~\eqref{eq:policy}) to estimate the infinite-horizon objective function.
See~\cite[Algorithm 5]{lan2022complexity} for more details.
Here, in-sample uses the same samples from the SAA problem~\eqref{eq:qgp_a3b} while out-of-sample uses new samples from the same distribution as~\eqref{eq:saa_xi}.

\begin{table}[h]
   \centering
   \begin{tabular}{@{}llrrrr@{}} \toprule
    Prob & Samples & CE-Inf-EDDP & EDDP & Cyc-SDDP \\ \midrule
    Inventory ($n=1$) & In & 3523.2 & \textbf{3523.1} & 3527.7 \\
    Inventory ($n=1$) & Out & \textbf{3538.4} & 3538.6 & 3542.1  \\ \midrule
    Hydro & In & \textbf{2.627e+8} & 3.022e+8 & 2.628e+08 \\
    Hydro & Out & 2.690e+8  & 3.211e+8 & \textbf{2.681e+8}  \\ \bottomrule
   \end{tabular}
   \caption{Evaluation is taken over $T'=1280$ periods with a discount factor $\lambda=0.9906$. Out-of-sample mean is shown over 30 seeds.
   The best (lowest) evaluations are bolded.
   } 
   \label{tab:hydro_final_oos}
\end{table}

The performances for both a 1-D inventory problem and a hydrothermal planning problem with $\lambda=0.9906$ are put in Table~\ref{tab:hydro_final_oos}. 
We see all methods have similar objectives for the inventory problem, which is expected since all methods can get nearly tight optimality gaps, as seen in Figure~\ref{fig:newsvendor_gaps}.
In the hydrothermal problem, our CE-Inf-EDDP attains matching in- and out-of-sample performance to Cyc-SDDP, which can be partially explained by the fact that both methods have similar lower bound estimates in Figure~\ref{fig:hydro_gaps}.
Moreover, EDDP has a worse (higher) objective, which aligns with its less tight lower bound from the same figure.
While both CE-Inf-EDDP and Cyc-SDDP have comparable performance, we should also recall that only the former can provably verify approximate optimality, as seen in Figure~\ref{fig:hydro_gaps}.

\section{Conclusion}
By leveraging the stationarity of the problem, our family of infinite-horizon dual dynamic programming methods matches state-of-the-art complexity and seems to exhibit promising numerical performance.
These \revtwo{advances} stem from the simplicity of a single cutting-plane model (via cut sharing), more frequent updates, and a new search point selection rule that encourages longer and adaptive exploration. 
\revone{An interesting future direction is to generalize these developments towards the infinite-horizon periodic case with periodicity $m$, where this paper considers $m=1$.
Since larger periods $m \gg 1$ result in fewer cuts shared between time periods, one promising question is to develop a relaxed cut sharing that can potentially share cuts between similar periods, such that those with similar costs, feasible sets, and/or random data.
}

\subsection*{Acknowledgments}
We thank Alex Shapiro for providing some comments on a preliminary draft that improved the clarity of this work.
We would also like to thank the reviewers for their comments and suggestions, and in particular, for an insightful comment about choosing a reset search point.

\bibliographystyle{plain}
\bibliography{biblio}

\end{document}

%% file: header.tex
\usepackage{graphicx}
\usepackage{amsmath}
\usepackage{hyperref} 
\usepackage{array}
\usepackage{tabularx,ragged2e,booktabs,caption}
\usepackage[shortlabels]{enumitem}
\usepackage{ulem}
\usepackage{subcaption}
\usepackage{float}
\usepackage{lipsum}
\usepackage{titling}
\usepackage{xcolor}

\usepackage{amsthm}
\usepackage{amssymb}

\usepackage{algorithm}
\usepackage{algpseudocode}

\usepackage{times}

\usepackage[
  a4paper,
  top=1in,
  bottom=1in,
  left=1.0in,
  right=1.0in,
]{geometry}
\usepackage{amsfonts}

\usepackage{caption}
\usepackage{subcaption}

\usepackage{bm}
\newcommand{\B}[1]{#1}

\theoremstyle{plain}
\newtheorem{theorem}{Theorem}[section]
\newtheorem{lemma}[theorem]{Lemma}
\newtheorem{corollary}[theorem]{Corollary}
\newtheorem{proposition}[theorem]{Proposition}

\theoremstyle{definition}
\newtheorem{definition}{Definition}[section]
\newtheorem{assumption}{Assumption}[section]
\theoremstyle{remark}

\newlength{\leftstackrelawd}
\newlength{\leftstackrelbwd}
\def\leftstackrel#1#2{\settowidth{\leftstackrelawd}%
{${{}^{#1}}$}\settowidth{\leftstackrelbwd}{$#2$}%
\addtolength{\leftstackrelawd}{-\leftstackrelbwd}%
\leavevmode\ifthenelse{\lengthtest{\leftstackrelawd>0pt}}%
{\kern-.5\leftstackrelawd}{}\mathrel{\mathop{#2}\limits^{#1}}}

\let\texdisplaystyle\displaystyle
\def\displaytotextstyle{\textstyle\let\displaystyle\texdisplaystyle}
\newenvironment{talign}
 {\let\displaystyle\displaytotextstyle\align}
 {\endalign}
\newenvironment{talign*}
 {\let\displaystyle\displaytotextstyle\csname align*\endcsname}
 {\endalign}

\usepackage[sort,numbers]{natbib}

\title{\textbf{Dual dynamic programming for stochastic
programs over an infinite horizon}\thanks{H. Milton Stewart School of Industrial and Systems Engineering,  Georgia Institute of Technology (Email: \texttt{calebju4@gatech.edu}, \texttt{george.lan@isye.gatech.edu})}}
\author{
Caleb Ju\thanks{Supported by Department of Energy Computational Science Graduate Fellowship under Award Number DE-SC0022158.}
\and
Guanghui Lan
}
\date{\vspace{-5ex}}

\newcommand{\tv}{\tilde{v}}
\newcommand{\hv}{\hat{v}}

\newcommand{\uV}{\underline{V}}
\newcommand{\oV}{\overline{V}}
\newcommand{\ov}{\overline{v}}
\newcommand{\uF}{\underline{F}}
\newcommand{\oF}{\overline{F}}

\newcommand{\xeps}{\bar x^k_i}

\newcommand{\Seps}{\mathcal{S}_\epsilon}
\newcommand{\mX}{\mathcal{X}}
\newcommand{\mXeps}{\mathcal{X}_\epsilon}
\newcommand{\Geps}{P_{\epsilon}}

\usepackage{tikz}
\newcommand{\tikzxmark}{%
\tikz[scale=0.23] {
    \draw[line width=0.7,line cap=round] (0,0) to [bend left=6] (1,1);
    \draw[line width=0.7,line cap=round] (0.2,0.95) to [bend right=3] (0.8,0.05);
}}
\newcommand{\tikzcmark}{%
\tikz[scale=0.23] {
    \draw[line width=0.7,line cap=round] (0.25,0) to [bend left=10] (1,1);
    \draw[line width=0.8,line cap=round] (0,0.35) to [bend right=1] (0.23,0);
}}

\providecommand{\keywords}[1]{ {\small \textbf{Keywords.} #1}}

\providecommand{\mscclass}[1]{ {\small \textbf{MSC subject classifications.} #1}}

\newcommand{\revone}[1]{#1} 
\newcommand{\revtwo}[1]{#1} 

%% file: biblio.bib
@article{pereira1991multi,
  title={Multi-stage stochastic optimization applied to energy planning},
  author={Pereira, Mario VF and Pinto, Leontina MVG},
  journal={Mathematical programming},
  volume={52},
  pages={359--375},
  year={1991},
  publisher={Springer}
}

@article{mcclain1977horizon,
  title={Horizon effects in aggregate production planning with seasonal demand},
  author={McClain, John O and Thomas, Joseph},
  journal={Management Science},
  volume={23},
  number={7},
  pages={728--736},
  year={1977},
  publisher={INFORMS}
}

@article{van2005comparison,
  title={A comparison of methods for lot-sizing in a rolling horizon environment},
  author={Van Den Heuvel, Wilco and Wagelmans, Albert PM},
  journal={Operations Research Letters},
  volume={33},
  number={5},
  pages={486--496},
  year={2005},
  publisher={Elsevier}
}

@article{jurek2011optimal,
  title={Optimal value and growth tilts in long-horizon portfolios},
  author={Jurek, Jakub W and Viceira, Luis M},
  journal={Review of Finance},
  volume={15},
  number={1},
  pages={29--74},
  year={2011},
  publisher={Oxford University Press}
}

@article{brown2011dynamic,
  title={Dynamic portfolio optimization with transaction costs: Heuristics and dual bounds},
  author={Brown, David B and Smith, James E},
  journal={Management Science},
  volume={57},
  number={10},
  pages={1752--1770},
  year={2011},
  publisher={INFORMS}
}

@article{cruise2019control,
  title={Control of energy storage with market impact: Lagrangian approach and horizons},
  author={Cruise, James and Flatley, Lisa and Gibbens, Richard and Zachary, Stan},
  journal={Operations Research},
  volume={67},
  number={1},
  pages={1--9},
  year={2019},
  publisher={INFORMS}
}

@article{lan2022complexity,
  title={Complexity of stochastic dual dynamic programming},
  author={Lan, Guanghui},
  journal={Mathematical Programming},
  volume={191},
  number={2},
  pages={717--754},
  year={2022},
  publisher={Springer}
}

@article{shapiro2011analysis,
  title={Analysis of stochastic dual dynamic programming method},
  author={Shapiro, Alexander},
  journal={European Journal of Operational Research},
  volume={209},
  number={1},
  pages={63--72},
  year={2011},
  publisher={Elsevier}
}

@article{philpott2008convergence,
  title={On the convergence of stochastic dual dynamic programming and related methods},
  author={Philpott, Andrew B and Guan, Ziming},
  journal={Operations Research Letters},
  volume={36},
  number={4},
  pages={450--455},
  year={2008},
  publisher={Elsevier}
}

@article{lan2024numerical,
  title={Numerical methods for convex multistage stochastic optimization},
  author={Lan, Guanghui and Shapiro, Alexander},
  journal={Foundations and Trends{\textregistered} in Optimization},
  volume={6},
  number={2},
  pages={63--144},
  year={2024},
  publisher={Now Publishers, Inc.}
}

@book{shapiro2021lectures,
  title={Lectures on stochastic programming: modeling and theory},
  author={Shapiro, Alexander and Dentcheva, Darinka and Ruszczynski, Andrzej},
  year={2021},
  publisher={SIAM},
  volume={3},
}

@book{nesterov2018lectures,
  title={Lectures on convex optimization},
  author={Nesterov, Yurii},
  year={2018},
  publisher={Springer},
  volume={137},
}

@article{dowson2020policy,
  title={The policy graph decomposition of multistage stochastic programming problems},
  author={Dowson, Oscar},
  journal={Networks},
  volume={76},
  number={1},
  pages={3--23},
  year={2020},
  publisher={Wiley Online Library}
}

@book{rockafellar1972convex,
  title={Convex Analysis},
  author={Ralph Tyrrell Rockafellar},
  year={1972},
  publisher={Princeton University Press}
}

@article{kleywegt2002sample,
  title={The sample average approximation method for stochastic discrete optimization},
  author={Kleywegt, Anton J and Shapiro, Alexander and Homem-de-Mello, Tito},
  journal={SIAM Journal on optimization},
  volume={12},
  number={2},
  pages={479--502},
  year={2002},
  publisher={SIAM}
}

@book{ben2001lectures,
  title={Lectures on modern convex optimization: analysis, algorithms, and engineering applications},
  author={Ben-Tal, Aharon and Nemirovski, Arkadi},
  year={2001},
  publisher={SIAM}
}

@article{zhang2020distributionally,
  title={On distributionally robust multistage convex optimization: New algorithms and complexity analysis},
  author={Zhang, Shixuan and Sun, Xu Andy},
  journal={arXiv preprint arXiv:2010.06759},
  year={2024}
}

@article{baucke2017deterministic,
  title={A deterministic algorithm for solving multistage stochastic programming problems},
  author={Baucke, Regan and Downward, Anthony and Zakeri, Golbon},
  journal={Optimization Online},
  pages={25},
  year={2017}
}

@article{zhang2022stochastic,
  title={Stochastic dual dynamic programming for multistage stochastic mixed-integer nonlinear optimization},
  author={Zhang, Shixuan and Sun, Xu Andy},
  journal={Mathematical Programming},
  volume={196},
  number={1},
  pages={935--985},
  year={2022},
  publisher={Springer}
}

@article{homem2011sampling,
  title={Sampling strategies and stopping criteria for stochastic dual dynamic programming: a case study in long-term hydrothermal scheduling},
  author={Homem-de-Mello, Tito and De Matos, Vitor L and Finardi, Erlon C},
  journal={Energy Systems},
  volume={2},
  number={1},
  pages={1--31},
  year={2011},
  publisher={Springer}
}

@article{girardeau2015convergence,
  title={On the convergence of decomposition methods for multistage stochastic convex programs},
  author={Girardeau, Pierre and Leclere, Vincent and Philpott, Andrew B},
  journal={Mathematics of Operations Research},
  volume={40},
  number={1},
  pages={130--145},
  year={2015},
  publisher={INFORMS}
}

@article{grinold1977finite,
  title={Finite horizon approximations of infinite horizon linear programs},
  author={Grinold, Richard C},
  journal={Mathematical Programming},
  volume={12},
  pages={1--17},
  year={1977},
  publisher={Springer}
}

@article{birge1985decomposition,
  title={Decomposition and partitioning methods for multistage stochastic linear programs},
  author={Birge, John R},
  journal={Operations research},
  volume={33},
  number={5},
  pages={989--1007},
  year={1985},
  publisher={INFORMS}
}

@article{shapiro2023dual,
  title={Dual bounds for periodical stochastic programs},
  author={Shapiro, Alexander and Cheng, Yi},
  journal={Operations Research},
  volume={71},
  number={1},
  pages={120--128},
  year={2023},
  publisher={INFORMS}
}

@article{shapiro2020periodical,
  title={Periodical multistage stochastic programs},
  author={Shapiro, Alexander and Ding, Lingquan},
  journal={SIAM Journal on Optimization},
  volume={30},
  number={3},
  pages={2083--2102},
  year={2020},
  publisher={SIAM}
}

@article{nannicini2021benders,
  title={A Benders squared ({$B^2$}) framework for infinite-horizon stochastic linear programs},
  author={Nannicini, Giacomo and Traversi, Emiliano and Calvo, Roberto Wolfler},
  journal={Mathematical Programming Computation},
  volume={13},
  number={4},
  pages={645--681},
  year={2021},
  publisher={Springer}
}

@article{guigues2020inexact,
  title={Inexact cuts in stochastic dual dynamic programming},
  author={Guigues, Vincent},
  journal={SIAM Journal on Optimization},
  volume={30},
  number={1},
  pages={407--438},
  year={2020},
  publisher={SIAM}
}

@book{powell2007approximate,
  title={Approximate Dynamic Programming: Solving the curses of dimensionality},
  author={Powell, Warren B},
  volume={703},
  year={2007},
  publisher={John Wiley \& Sons}
}

@book{caines2018linear,
  title={Linear stochastic systems},
  author={Caines, Peter E},
  year={2018},
  publisher={SIAM}
}

@article{chen1999convergent,
  title={Convergent cutting-plane and partial-sampling algorithm for multistage stochastic linear programs with recourse},
  author={Chen, Zhi-Long and Powell, Warren B},
  journal={Journal of Optimization Theory and Applications},
  volume={102},
  number={3},
  pages={497--524},
  year={1999},
  publisher={Springer}
}
